\documentclass[11pt]{article}

\usepackage{amsmath}
\usepackage{amssymb}
\usepackage{amscd}

\usepackage[english]{babel}

\usepackage[raiselinks=false,plainpages=false,pdfpagelabels,naturalnames=true,colorlinks=true,citecolor=blue,urlcolor=blue,linkcolor=blue,bookmarksopen=true,pagebackref,hyperindex]{hyperref}

\allowdisplaybreaks[1]

\usepackage{latexsym}
\usepackage[emu]{cmll}
\usepackage{prooftree}

\def\Hide#1{\relax}

\DeclareSymbolFont{AMSb}{U}{msb}{m}{n}
\DeclareSymbolFontAlphabet{\mathbb}{AMSb}
\DeclareSymbolFont{symbolsC}{U}{txsyc}{m}{n}
\DeclareMathSymbol{\rJoin}{\mathrel}{symbolsC}{89}

\def\BB{\mathbb{B}}

\def\NN{\mathbb{N}}

\def\QQ{\mathbb{Q}}
\def\RR{\mathbb{R}}

\def\cB{\mathcal{B}}
\def\cC{\mathcal{C}}

\def\cI{\mathcal{I}}

\def\VA{\mathbf{A}}

\def\VC{\mathbf{C}}
\def\VD{\mathbf{D}}
\def\VE{\mathbf{E}}
\def\VF{\mathbf{F}}
\def\VG{\mathbf{G}}
\def\VH{\mathbf{H}}

\def\VK{\mathbf{K}}
\def\VL{\mathbf{L}}
\def\VM{\mathbf{M}}

\def\VP{\mathbf{P}}
\def\VQ{\mathbf{Q}}
\def\VR{\mathbf{R}}
\def\VS{\mathbf{S}}
\def\VT{\mathbf{T}}
\def\VU{\mathbf{U}}

\def\VY{\mathbf{Y}}

\def\Func#1{{\sf{#1}}}

\def\And{\land}
\def\iAnd{\otimes} 

\def\Con#1{\Func{Con}\,#1}

\def\Diff{\mathop{\backslash}}

\def\Id{\Func{id}}
\def\Im{\Func{im}}
\def\Imp{\Rightarrow}

\def\Ker{\Func{ker}}

\def\Max{\Func{max}}
\def\Min{\Func{min}}

\def\Not{\lnot}
\def\Or{\lor}

\def\ST{\mathrel{|}}

\def\Th{\Func{Th}}
\def\To{\rightarrow}

\def\all#1{\forall #1 {\cdot}\:}

\def\UI{\mathbf{[0, 1]}}
\def\RRnn{\RR^{{\ge}0}}

\newtheorem{Theorem}{Theorem}[subsection]
\newtheorem{Lemma}[Theorem]{Lemma}

\newtheorem{Definition}{Definition}[subsection]
\newtheorem{Remark}{Remark}[subsection]
\newtheorem{Example}{Example}[subsection]
\def\Proof{\par \noindent{\bf Proof: }}
\def\Done{\hfill\rule{0.5em}{0.5em}}
\def\Var{\Func{Var}}

\def\dotplus{\mathop{\dot{+}}}

\def\lL{{\cal L}}

\def\MONE{[{\sf m}_1]}
\def\MTWO{[{\sf m}_2]}
\def\MTHREE{[{\sf m}_3]}
\def\Mi{[{\sf m}_i]}
\def\OONE{[{\sf o}_1]}
\def\OTWO{[{\sf o}_2]}
\def\OTHREE{[{\sf o}_3]}
\def\OFOUR{[{\sf o}_4]}
\def\Oj{[{\sf o}_j]}
\def\BOUNDEDBELOW{[{\sf b}]}
\def\RESIDUAL{[{\sf r}]}
\def\ONEAX{[{\sf ann}]}
\def\INV{[{\sf dne}]}
\def\HOOP{[{\sf cwc}]}

\def\ClCl#1{\mathsf{IR}(#1)}

\def\BPoc{{\Func{Poc}_1}}

\def\Set{\Func{Set}}
\def\Hom{\Func{Hom}}
\def\Ass{\Func{Ass}}
\def\Sem{\Func{Sem}}

\def\Lnot{{{}^{\perp}}}
\def\Lolly{\multimap}

\def\axWk{(\Func{Wk})}
\def\axComp{(\Func{Comp})}
\def\axComm{(\Func{Comm})}
\def\axUncurry{(\Func{Uncurry})}
\def\axCurry{(\Func{Curry})}
\def\axEFQ{(\Func{EFQ})}
\def\axDNE{(\Func{DNE})}
\def\axCWC{(\Func{CWC})}
\def\axCon{(\Func{Con})}

\def\Logic#1#2{\mbox{{\bf #1}}_{\mbox{\bf #2}}}

\def\ALi{\Logic{AL}{i}}
\def\ALm{\Logic{AL}{m}}
\def\BL{\Logic{BL}{}}
\def\ALc{\Logic{AL}{c}}
\def\Luk{\Logic{{\L}}{\relax}}
\def\LLi{\Logic{{\L}L}{i}}
\def\LLm{\Logic{{\L}L}{m}}
\def\LLc{\Logic{{\L}L}{c}}

\def\IL{\Logic{IL}{}}
\def\ML{\Logic{ML}{}}

\def\ordSum{\mathop{\stackrel{\frown}{\relax}}}
\def\rImp{\mathop{\rightarrow}}

\newcommand{\kTrans}[1]{{#1}^{\sf K}}

\newcommand{\ggTrans}[1]{{#1}^{\sf Gen}}
\newcommand{\glTrans}[1]{{#1}^{\sf Gli}}

\newcommand{\kSem}{\kTrans{\mu}}
\newcommand{\ggSem}{\ggTrans{\mu}}
\newcommand{\glSem}{\glTrans{\mu}}
\newcommand{\sSem}{\mu^{\sf S}}

\title{On Pocrims and Hoops}

\author{Rob Arthan \& Paulo Oliva}

\begin{document}
\maketitle

\begin{abstract}
Pocrims and suitable specialisations thereof are structures that
provide the natural algebraic semantics for a minimal affine logic and its
extensions.  Hoops comprise a special class of pocrims that provide
algebraic semantics for what we view as an intuitionistic analogue of the
classical multi-valued {\L}ukasiewicz logic.  We present some contributions to
the theory of these algebraic structures.  We give a new
proof that the class of hoops is a variety.  
We use a new indirect method to establish several important identities in the
theory of hoops: in particular, we prove that the double negation mapping in a
hoop is a homormorphism.  This leads to an investigation of algebraic analogues
of the various double negation translations that are well-known from proof
theory.  We give an algebraic framework for studying the semantics of double
negation translations and use it to prove new results about
the applicability of the double negation translations due to Gentzen and
Glivenko.

\end{abstract}

\section{Introduction}\label{sec:introduction}

Pocrims provide the natural algebraic models for a minimal affine logic,
$\ALm$, while hoops provide the models for what we view as a minimal
analogue, $\LLm$, of {\L}ukasiewicz's classical infinite-valued logic
$\LLc$.  This paper presents some new results on the algebraic structure
of pocrims and hoops.  Our main motivation for this work is in the
logical aspects: we are interested in criteria for provability in
$\ALm$, $\LLm$ and related logics. We  develop a useful practical test
for provability in $\LLm$ and apply it to a range of problems including
a study of the various double negation translations in these logics.

We begin in Section~\ref{sec:background} with a brief
introduction to the logical background and then give
the definitions and basic theory of the algebraic structures.
Since we believe the algebraic approach will be unfamiliar to some
readers who share our interest in the logical issues, this part
of the paper is largely expository, bringing together material that
is scattered over the literature.
We illustrate the material with a number of examples, including
all pocrims of order 2, 3 and 4. Some of our later results depend
on the existence of finite pocrims satisfying or failing to
satisfy certain identities: the witnesses were all found using
the late Bill McCune's Mace4 program~\cite{prover9-mace4},
which has proved an invaluable tool in our work.

In Section~\ref{sec:alg-sem} we review the algebraic semantics for the logics
introduced in Section~\ref{sec:background} and prove the soundness and
completeness of pocrims and appropriate specialisations thereof to these
logics. Again this section is largely expository however it concludes,
with a new proof that the class of hoops is a variety. The proof
provides an algorithm for translating a proof tree in the logic $\LLm$
into a semantically equivalent equational proof.

The equational theory of hoops is known to be decidable and it follows
from work of Bova and Montagna \cite{Bova:2009} that the decision
problem is in PSPACE. Unfortunately, their decision procedure is
infeasible in practice, even on small examples.  In
Section~\ref{sec:indirect-method}, we attempt to mitigate this
difficulty.  We begin by reviewing known results on the equational
decision problem for involutive hoops (i.e., hoops that satisfy an
algebraic analogue of the law of double negation elimination).  The
variety of involutive hoops can be shown to be definitionally equivalent
to the well-known variety of MV-algebras and the equational theory of
MV-algebras reduces to the theory of linear real arithmetic.  We then
reduce the decision problem for an identity in a general hoop to
restricted classes of finitely generated hoops enjoying special
algebraic properties.  This falls short of a decision procedure, but
provides an efficient heuristic that can be used to prove many important
identities, whose formal proofs, if known, are extremely intricate.  We
give several interesting applications of this method, e.g., we show that
the set of idempotent elements in a hoop is the universe of a subhoop.

In Section~\ref{sec:neg-trans}, we use the method of
Section~\ref{sec:indirect-method} to show that the double negation mapping in a
hoop is a homomorphism.  We undertake an algebraic investigation of the double
negation translations of Kolmogorov, Gentzen and Glivenko.  Kolmogorov's
translation is shown to be correct for any extension of affine logic. The
Gentzen and Glivenko translations are correct for intuitionistic {\L}ukasiewicz
logic, but there are weaker extensions of affine logic for which Gentzen is
correct but Glivenko is not and {\it vice versa}.

%

\section{Background}\label{sec:background}

While the main emphasis of this paper is on algebraic structures,
our main motivation for studying those structures stems from an
interest in certain substructural propositional logics. We now define those logics.

\subsection{Nine Logics}\label{sec:nine-logics}

We work in a language, $\lL$, built from a countable set of variables
$\Var = \{V_1, V_2, \ldots\}$, the constant $1$ (falsehood) and the binary
connectives $\Lolly$ (implication) and $\iAnd$ (conjunction).
We write $A\Lnot$ for $A \Lolly 1$ and $0$ for $1 \Lolly 1$.
Our choice of notation for connectives is that commonly used for
affine logic, since all the systems
we consider are extensions of intuitionistic affine logic.
Our use of $1$ rather than $0$ for falsehood is taken from continuous logic
\cite{ben-yaacov-pedersen09}, which motivated our work in this area.

\begin{figure}[t]
\[
\begin{array}{|c|l|} \hline
\axComp& (A \Lolly B) \Lolly (B \Lolly C) \Lolly (A \Lolly C) \\[1mm] \hline
\axComm& A \iAnd B \Lolly B \iAnd A \\[1mm] \hline
\axCurry& (A \iAnd  B \Lolly C) \Lolly (A \Lolly B \Lolly C) \\[1mm] \hline
\axUncurry& (A \Lolly B \Lolly C) \Lolly  (A \iAnd B \Lolly C) \\ \hline
\axWk& A \iAnd B \Lolly A \\[1mm] \hline
\axEFQ& 1 \Lolly A \\[1mm] \hline
\axDNE& A \Lnot\Lnot \Lolly A \\[1mm] \hline
\axCWC& A \iAnd (A \Lolly B) \Lolly B \iAnd (B \Lolly A) \\[1mm] \hline
\axCon& A \Lolly A \iAnd A \\[1mm] \hline
\end{array}
\]
\caption{Axiom Schemata}
\label{fig:axioms}
\end{figure}

As usual, we adopt the convention that $\Lolly$ associates to the right and has
lower precedence than $\iAnd$, which in turn has lower precedence than $(\cdot)\Lnot$.
So, for example, the brackets in $(A \iAnd  (B\Lnot)) \Lolly (C \Lolly (D \iAnd F))$ are all redundant, while those in $(((A \Lolly B) \Lolly C) \iAnd D)\Lnot$ are all required.

If $T$ is a set of formulas in the language $\lL$, the {\em deductive closure},
$\overline{T}$,
of $T$ is the smallest subset of $\lL$ that contains $T$ and is closed
under {\em modus ponens} (i.e., if $A \in \overline{T}$ and $A \Lolly B \in \overline{T}$ then $B \in \overline{T}$).
If $T = \overline{T}$, we say $T$ is {\em deductively closed} or {\em a theory}.
For our purposes in this paper a logic is just a theory.
However, we will often write ``$T$ proves $A$'' or ``$A$ is derivable in $T$'' as
a suggestive alternative to $A \in T$.
If $S$ and $T$ are sets of formulas, e.g., theories or axiom schemata,
we write $S + T$ for $\overline{S \cup T}$.

We will consider nine axiom schemata as shown in the table of Figure \ref{fig:axioms}. These are:
{\em composition, commutativity of conjunction,
currying, uncurrying, weakening, {\it ex falso quodlibet},
double negation elimination, commutativity of weak
conjunction and contraction}.

We then consider nine combinations of these axiom schema, as shown in Figure \ref{fig:logics}. $\ALm$, $\ALi$, $\ALc$, $\LLm$, $\LLi$ and $\LLc$ are
minimal, intuitionistic and classical variants
of affine logic and {\L}ukasiewicz logic.
$\ML$, $\IL$ and $\BL$ have both weakening,
$\axWk$, and contraction, $\axCon$,
and so are the implication-conjunction fragments of the usual minimal,
intuitionistic and boolean logics.
Over $\ALm$, the schema $\axCon$ implies the schema $\axCWC$.
In fact, as discussed in~\cite{arthan-oliva14a},
one can interpret $\axCWC$ as a weak form of the contraction rule.
We can consequently depict our nine logics in
the 2-dimensional diagram shown in Figure~\ref{fig:logics}
(in which the rectangles are push-outs in the poset of
deductively closed subsets of $\lL$).

\begin{figure}[t]
\[
\begin{array}{|c|l|} \hline
\ALm & \axComp + \axComm + \axCurry + \axUncurry + \axWk \\[1mm] \hline
\ALi & \ALm + \axEFQ \\[1mm] \hline
\ALc & \ALi + \axDNE \\[1mm] \hline
\LLm & \ALm + \axCWC \\[1mm] \hline
\LLi & \LLm + \axEFQ \\[1mm] \hline
\LLc & \LLi + \axDNE \\[1mm] \hline
\ML  & \ALm + \axCon \\[1mm] \hline
\IL  & \ML  + \axEFQ \\[1mm] \hline
\BL  & \IL  + \axDNE \\[1mm] \hline
\end{array}
\]
\caption{Logics}
\label{fig:logics}
\end{figure}

\begin{figure}[h]
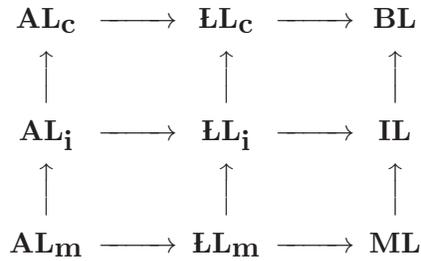

\[
\begin{CD}
\ALc  @>>> \LLc  @>>> \BL \\
@AAA       @AAA       @AAA \\
\ALi  @>>> \LLi  @>>> \IL \\
@AAA       @AAA       @AAA \\
\ALm  @>>> \LLm  @>>> \ML \\
\end{CD}
\]
\caption{Relationships between the Logics}
\label{fig:logics}
\end{figure}

It was shown in the 1950s by Rose and Rosser \cite{rose-rosser58} (and also using a different method of proof by Chang \cite{Chang59}) that the Hilbert-style system
$\Luk$ with the following axiom schemata\footnote{
Following {\L}ukasiewicz, Rose and Rosser used Polish notation.
Rose and Rosser write $CAB$ for our $A\Lolly B$ and $NA$ for our $A\Lnot$.
Chang followed this in the relatively few fragments of syntax that appear
in his treatment.
} is sound and complete for {\L}ukasiewicz's many-valued logical system
where the truth values are real numbers in the interval $[0, 1]$
and where $\Lolly$ and $(\cdot)\Lnot$ are modelled by $(a, b) \mapsto \Min(a + b, 1)$ and $a \mapsto 1 - a$ respectively\footnote{
Throughout this paper we adopt the convention that truth values are ordered
by increasing logical strength, so $0$ represents truth and $1$ represents
falsehood.
}.
\[
\begin{array}{rl}
\mbox{(A1)} & A \Lolly (B \Lolly A) \\
\mbox{(A2)} & (A \Lolly B) \Lolly (B \Lolly C) \Lolly (A \Lolly C) \\
\mbox{(A3)} & ((A \Lolly B) \Lolly B) \Lolly ((B \Lolly A) \Lolly A) \\
\mbox{(A4)} & (A\Lnot \Lolly B\Lnot) \Lolly (B \Lolly A) \\
\end{array}
\]
Note that in $\Luk$ conjunction can be defined as $A \iAnd B \mathrel{:\equiv} (A\Lnot \Lolly B)\Lnot$.
We will see in Section~\ref{sec:id-inv-hoops} that our
$\LLc$ is equivalent to $\Luk$.

\subsection{Pocrims and Hoops}\label{sec:the-algebras}

\begin{Definition} A {\em pocrim}\footnote{The name is an acronym for
``partially ordered, commutative, residuated, integral monoid'',
Strictly speaking, this is a {\em dual} pocrim,
since we order it by increasing logical strength and write it additively.} $\VP$
is a structure for the signature $(0, +, \rImp)$ of type
$(0, 2, 2)$ satisfying the following laws, in which $x \ge y$
is an abbreviation for $x \rImp y = 0$:
\begin{align*}
&(x + y) + z = x + (y + z) \tag*{$\MONE$} \\
&x + y = y + x \tag*{$\MTWO$} \\
&x + 0 = x \tag*{$\MTHREE$} \\
&x \ge x \tag*{$\OONE$} \\
&\mbox{if $x \ge y$ and $y \ge z$, then $x \ge z$} \tag*{$\OTWO$} \\
&\mbox{if $x \ge y$ and $y \ge x$, then $x = y$} \tag*{$\OTHREE$} \\
&\mbox{if $x \ge y$, then $x + z \ge y + z$} \tag*{$\OFOUR$} \\
&x \ge 0 \tag*{\BOUNDEDBELOW} \\
&\mbox{$x + y \ge z$ iff $x \ge y \rImp z$.} \tag*{\RESIDUAL}
\end{align*}
\end{Definition}
We will see that pocrims provide models for our logics:
$\rImp$ is the semantic counterpart of the syntactic implication $\Lolly$, whereas $+$ corresponds to the syntactic conjunction $\iAnd$.
As with the syntactic connectives,
we adopt the convention that $\rImp$ associates
to the right and has lower precedence than $+$.
So the brackets in $x + (x \rImp y)$ are necessary while those
in $x \rImp (y \rImp z)$ may be omitted.
Throughout this paper, we adopt the convention that if $\VP$ is a structure then $P$ is its universe.

The laws~$\Mi$, $\Oj$ and $\BOUNDEDBELOW$ say that
$(P; 0, +; {\ge})$ is a partially ordered commutative monoid
with the identity $0$ as least element.
Law~$\RESIDUAL$, the {\em residuation property},
says that for any $y$ and $z$ the
set $\{x \ST x + y \ge z\}$ is non-empty and has $y \rImp z$
as least element.
Taking $x = y \rImp z$ in $\RESIDUAL$ and using $\OONE$, we have
that $(y \rImp z) + y \ge z$, an algebraic analogue of {\em modus ponens}.

A pocrim is said to be {\em bounded} if it has a (necessarily unique)
{\em annihilator}, i.e., an element $1$ such that for every $x$ we have:
\begin{align*}
& 1 = x + 1. \tag*{\ONEAX}
\end{align*}
In a bounded pocrim $\VP$, we have that $1 = x + 1 \ge x + 0 = x$ for any $x$,
so that $(M; \ge)$ is indeed a bounded ordered set.
We write $\Not x$ for $x \rImp 1$ (and give  $\Not$ higher precedence than the binary operators).
Note that any finite pocrim $\VP$ is bounded, the annihilator being given by
$\sum_{x \in P} x$.

A pocrim is said to be {\em involutive} if it is bounded and satisfies
the double-negation identity:
\begin{align*}
&\Not\Not x = x. \tag*{\INV}
\end{align*}
We will often write $\delta(x)$ for $\Not\Not x$.
In any bounded pocrim, the set $\{0, 1\}$ is closed under $+$ and $\rImp$ and
so, as $\Not 0 = 1$ and $\Not 1 = 0$, $\{0, 1\}$ is the universe
of an involutive subpocrim.

\begin{Example}
There is a unique pocrim $\BB$ with two elements.
It is involutive and provides the standard model for classical Boolean logic.
\end{Example}

If $x$ and $y$ are elements of a pocrim, $x + (x \rImp y)$ is an upper bound for $x$ and $y$ as is $y + (y \rImp x)$.
Logically, we can view either of these as a weak form of conjunction.
Pocrims in which the two upper bounds coincide turn out to have many pleasant
properties, motivating the following definition.

\begin{Definition}[B\"{u}chi \& Owens\cite{BO}]
A {\em hoop}\footnote%
{
B\"{u}chi and Owens \cite{BO} write of hoops that ``their importance \ldots
merits recognition with a more euphonious name than the merely descriptive
``commutative complemented monoid''''. Presumably they chose ``hoop'' as a
euphonious companion to ``group'' and ``loop''.
}
is a pocrim that satisfies {\em commutativity of weak conjunction}:
\begin{align*}
&x + (x \rImp y) = y + (y \rImp x). \tag*{\HOOP}
\end{align*}
\end{Definition}

The following lemma provides some useful characterisations of hoops.

\begin{Lemma}
If $\VP$ is a pocrim, the following are equivalent:
\begin{enumerate}

\item \label{hoop-cond}
$\VP$ is a hoop. I.e., $\VP$ satisfies $x + (x \rImp y) = y + (y \rImp x)$.
\item \label{nat-order}
$\VP$ is {\em naturally ordered}. I.e., 
for every $x, y \in P$ such that $x \ge y$, there is $z \in P$
such that $x = y + z$.
\item \label{hoop-cond-cond}
For every $x, y \in P$ such that $x \ge y$, $x = y + (y \rImp x)$.
\item \label{hoop-cond-ge}
$\VP$ satisifies $x + (x \rImp y) \ge y + (y \rImp x)$
\end{enumerate}
\end{Lemma}
\Proof
\emph{\ref{hoop-cond}} $\Imp$ \emph{\ref{nat-order}}:
Assume that $\VP$ satisfies $x + (x \rImp y) = y + (y \rImp x)$
and that $x, y \in P$ satisfy $x \ge y$, i.e., $x \rImp y = 0$.
Taking $z = y \rImp x$, we have
$x = x + 0 = x + (x \rImp y) = y + (y \rImp x) = y + z$.

\noindent
\emph{\ref{nat-order}} $\Imp$ \emph{\ref{hoop-cond-cond}}:
Assume that $\VP$ is naturally ordered
and that $x, y \in P$ satisfy $x \ge y$.
Then $x = y + z$ for some $z$.
By the residuation property, we have $z \ge y \rImp x$,
hence $x = y + z \ge y + (y \rImp x) \ge x$ and
so $x = y + (y \rImp x)$.

\noindent
\emph{\ref{hoop-cond-cond}} $\Imp$ \emph{\ref{hoop-cond-ge}}:
assume that $\VP$ satisfies $x = y + (y \rImp x)$ whenever
$x, y \in P$ and $x \ge y$. Given any $x, y \in P$,
we have $x + (x \rImp y) \ge y$, whence
$x + (x \rImp y) = y + (y \rImp x + (x \rImp y)) \ge y + (y \rImp x)$.

\noindent
\emph{\ref{hoop-cond-ge}} $\Imp$ \emph{\ref{hoop-cond}}:
exchange $x$ and $y$ and use the fact that $\ge$ is antisymmetric.
\Done \\

We will now give an outline of some basic algebraic properties of pocrims and
hoops omitting most of the proofs.  See~\cite{raftery07} for further
information about pocrims in general and involutive pocrims in particular.
See~\cite{blok-ferreirim00} for further information about hoops.

Given a linearly ordered abelian group $\VG$, there is a hoop
$\VG^{{\ge}0} = (\{x : G \ST x \ge 0\}; 0, +, \rImp)$, where
$x \rImp y = \Max(0, y - x)$. So for example, taking
$G$ to be the additive group of real numbers, we have the
hoop $\RR^{{\ge}0}$ whose elements are non-negative real numbers.
Given an element $a$ of a linearly ordered hoop $\VH$, there is a bounded hoop
$\VH^{\leq a} = (\{x : H \ST x \le a\}; 0, +_{a}, \rImp)$
where $x \mathop{+_{a}} y = \Min(a, x + y)$ so that $a$ becomes the annihilator.
If we compose these two constructions, the resulting bounded hoop
$\VG^{[0, a]}$ is involutive, since it satisfies $\Not x = a - x$.

\begin{Example}\label{eg:ui}
We write $\UI$ for the involutive hoop $\RR^{[0, 1]}$ obtained by the above
constructions taking $\VG = (\RR; 0, +)$ and $a = 1$.
Thus the universe of $\UI$ is the unit interval and the operations are
given by:
\[ x \dotplus y = \Min(x + y, 1) \quad \quad \quad x \rImp y = \Max(y - x, 0) \]
(where we write $\dotplus$ rather than $+$ for the hoop operation to
distinguish it from addition of real numbers).
$\UI$ provides an infinite model of classical {\L}ukasiewicz logic $\LLc$
(as does $\VG^{[0, 1]}$ for any dense subgroup of $(\RR; 0, +)$
containing 1).
\end{Example}

\begin{Example}\label{eg:vl}
For any integer $m \ge 1$, define $\VR_m$ to be the additive subgroup of $\QQ$
generated by $\frac{1}{m}$.
For $n \ge 2$, let $\VL_n = \VR_{n-1}^{[0, 1]}$.
Thus the universe of $\VL_n$ is $L_n = \{0, \frac{1}{n-1}, \frac{2}{n-1}, \ldots, \frac{n-2}{n-1}, 1\}$ and
the operations $+$ and $\rImp$ on $L_n$ are given by the same formulas
as for $\UI$ in Example~\ref{eg:ui}.
The hoops $\VL_n$ are involutive and provide natural finite
models of classical {\L}ukasiewicz logic $\LLc$.
\end{Example}

A hoop $\VH$ is said to be a {\em Wajsberg} hoop if it satisfies $(x \rImp y)
\rImp y = (y \rImp x) \rImp x$).  This is the algebraic equivalent of the axiom
schema (A3) of Section~\ref{sec:nine-logics}.  A bounded hoop is Wajsberg iff
it is involutive. There are, however, unbounded Wajsberg hoops, for instance:

\begin{Example}\label{eg:rrnn}
The unbounded hoop $\RRnn$ is Wajsberg. In fact, in $\RR^{{\ge}0}$
$(x \rImp y) \rImp y$ and $(y \rImp x) \rImp x$ are both equal to
$\Min(x, y)$.
\end{Example}

If $\VC$ and $\VD$ are pocrims, the {\em ordinal sum},
$\VC \ordSum \VD$, is the pocrim $(C \sqcup (D \Diff \{0\}), 0, +, \rImp)$
where $+$ and $\rImp$ extend the given operations on $C$ and $D$
to the disjoint union $C \sqcup (D \Diff \{0\})$ in such a way that
whenever $c \in C$ and $0 \not= d \in D$, $c + d = d$ (implying that
$c \rImp d = d$ and $d \rImp c = 0$).
Thus the order type of $\VC \ordSum \VD$ is the concatenation of the
partial orders $(C; \ge)$ and $(D \Diff \{0\}; \ge)$.
If $D \not= \{0\}$,
$\VC \ordSum \VD$ is bounded iff $\VD$ is bounded
and can only be involutive if $C = \{0\}$, since if $0 \not= c \in C$,
then, in $\VC \ordSum \VD$, we have $\Not c = 1$, so that $\Not\Not c = 0 \not= c$.
$\VC \ordSum \VD$ is a hoop iff both $\VC$ and $\VD$ are hoops.

\begin{Example}
Apart from $\VL_3$ there is one other pocrim with 3 elements,
namely $\VG_3 = \BB \ordSum \BB$.
$\VG_3$ is the first non-Boolean example in the sequence of idempotent pocrims
defined by the equations $\VG_2 = \BB$ and $\VG_{n+1} = \VG_n \ordSum \BB$.
$G_n = \{0, x_1, x_2, \ldots, x_{n-2}, 1\}$ with
$0 < x_1 < x_2 \ldots < x_{n-2} < 1$ and with operations defined by
\[
x + y = \Max\{x, y\} \quad\quad\quad
x \rImp y = \left\{
   \begin{array}{l@{\quad\quad}l}
      y & \mbox{if $y > x$} \\
      0 & \mbox{otherwise}
   \end{array} \right.
\]
The $\VG_n$ are finite models of intuitionistic propositional logic $\IL$.
They were used by G\"{o}del to prove that $\IL$
requires infinitely many truth values \cite{goedel32}.
In $\VG_n$, $\Not x = 1$ unless $x = 1$, so for $n > 2$, $\VG_n$ is not
involutive.
\end{Example}

\begin{Example}\label{eg:p-four}\label{eg:q-four}
It can be shown that there are 7 pocrims with 4 elements:
$\BB \times \BB$, $\VL_4$, $\VG_4$, $\BB \ordSum \VL_3$,
$\VL_3 \ordSum \BB$, $\VP_4$ and $\VQ_4$. $\VP_4$ and
$\VQ_4$ are the smallest pocrims that are not hoops
and are as follows: \\
$\VP_4$ comprises the chain $0 < p < q < 1$.
The operation tables for $\VP_4$ are as follows.
\[
\begin{array}{l@{\quad\quad}l@{\quad\quad}l}
\begin{array}{c|cccc}
   {+} & 0 & p & q & 1 \\\hline
    0  & 0 & p & q & 1 \\
    p  & p & 1 & 1 & 1 \\
    q  & q & 1 & 1 & 1 \\
    1  & 1 & 1 & 1 & 1
\end{array}
&
\begin{array}{c|cccc}
   {\rImp} & 0 & p & q & 1 \\\hline
    0      & 0 & p & q & 1 \\
    p      & 0 & 0 & p & p \\
    q      & 0 & 0 & 0 & p \\
    1      & 0 & 0 & 0 & 0
\end{array}
&
\begin{array}{c|c}
   \multicolumn{2}{c}{\delta}  \\
   \hline
    0      & 0 \\
    p      & p \\
    q      & p \\
    1      & 1
\end{array}
\end{array}
\]
(where for future reference we also tabulate the double negation mapping, $\delta$).
In $\VP_4$, $\delta(q) = p$, so $\VP_4$ is not involutive.
Moreover $\VP_4$ is not a hoop since it is not naturally ordered:
there is no $z$ with $p + z = q$.
However, the image of double negation is a subpocrim with
universe $\{0, p, 1\}$ isomorphic to the involutive
hoop $\VL_3$. \\
$\VQ_4$ comprises the chain $0 < u < v < 1$ and
has operation tables as follows:
\[
\begin{array}{l@{\quad\quad}l@{\quad\quad}l}
\begin{array}{c|cccc}
   {+} & 0 & u & v & 1 \\\hline
    0  & 0 & u & v & 1 \\ 
    u  & u & u & 1 & 1 \\ 
    v  & v & 1 & 1 & 1 \\ 
    1  & 1 & 1 & 1 & 1\end{array}
&
\begin{array}{c|cccc}
   {\rImp} & 0 & u & v & 1 \\\hline
    0      & 0 & u & v & 1 \\
    u      & 0 & 0 & v & v \\ 
    v      & 0 & 0 & 0 & u \\
    1      & 0 & 0 & 0 & 0
\end{array}
&
\begin{array}{c|c}
   \multicolumn{2}{c}{\delta} \\
   \hline
    0      & 0 \\
    u      & u \\
    v      & v \\
    1      & 1 
\end{array}
\end{array}
\] 
Like $\VP_4$, $\VQ_4$ is not naturally ordered
and hence not a hoop, because there is no $z$ with $u + z = v$.
$\VQ_4$ is involutive.
\end{Example}

An {\em ideal} in a hoop $\VH$ is a subset that forms the universe
of a downwards closed subhoop.
Trivially $H$ itself and $\{0\}$ are ideals.
If $X \subseteq H$, the ideal generated by $X$ comprises the set
of all $y \in H$ such that $y \le x_1 + \ldots + x_k$ for some
list $x_1$, \ldots, $x_k$ of elements of $X$.

If $f: \VH \To \VK$ is a hoop homomorphism, we define $\Ker(f)$, the
{\em kernel} of $f$, by $\Ker(f) = \{x : H \ST f(x) = 0\}$.
It is easy to verify that $\Ker(f)$ is an ideal.
Conversely, given an ideal $I$ in $\VH$,
the relation $\theta$ on $H$ defined by $x \mathrel{\theta} y$
iff $(x \rImp y) + (y \rImp x) \in I$ defines a congruence on $\VH$
such that if $p : \VH \To \VH/\theta$ is the natural projection of $\VH$
onto the quotient hoop\footnote{
We shall show in Section~\ref{thm:luk-provability-equational}
that the class of hoops is a variety,
so the quotient of a hoop by a congruence is in fact a hoop.
} $\VH/\theta$, then $\Ker(p) = I$.
This gives an isomorphism between the lattice of congruences on $\VH$
and its lattice of ideals.
We write $\VH/I$ for the quotient of $\VH$ by the congruence corresponding
to the ideal $I$.

\begin{Example}
If $\VC$ and $\VD$ are hoops, $C$ is an ideal in the ordinal sum $\VC \ordSum \VD$
and the quotient $(\VC \ordSum \VD)/C$ is isomorphic to $\VD$ via
an isomorphism that is left inverse to the natural inclusion of $\VD$
in $\VC \ordSum \VD$.
\end{Example}

If $\VP$ is a pocrim, $n \in \NN$ and $x \in P$, we
write $nx$ for the sum $\sum_{i=1}^n x$. $\VP$
is said to be {\em archimidean} if whenever $x, y \in \VP \Diff \{0\}$,
there is $n \in \NN$ such that $nx \ge y$.
By the equivalence between ideals and congruences, a hoop $\VH$ is simple,
i.e., admits no non-trivial congruences, iff the ideal generated by
any non-zero element of $H$ is $H$ itself.
It follows that a hoop is simple iff it is archimedean.
So, for example, the $\VL_n$ are all simple, while $\VG_n$ is simple
iff $n = 2$.

\begin{Example}
Let $\VE$ be the plane $\RR \times \RR$ given the structure of
a linearly ordered abelian group under vector addition and the
lexicographic ordering and let $a = (1, 0)$.
$\VE^{[0,a]}$ is not archimedean: the elements of $E^{[0,a]}$
on the $y$-axis form
a subhoop, $\VY$, such that $ny < a$ for any $n \in \NN$ and $y \in Y$.
The projection $\pi_1 : E \To \RR$ onto the $x$-axis induces
a surjective hoop homomorphism $f : \VE^{[0, a]} \To \UI$ and $\Ker(f) = Y$.
$\UI$ and $\VY$ are clearly archimedean, and hence simple.
It follows that $Y$ is the only non-trivial ideal in $\VE^{[0,a]}$.
\end{Example}

Recall that an algebra $\VA$ is said to be {\em subdirectly irreducible} if the
intersection $\Psi = \bigcap(\Con\;\VA \Diff \Delta)$ of all its congruences
other than the identity congruence, $\Delta$, is not equal to $\Delta$.  Using
the correspondence between congruences and ideals, $\VE^{[0,a]}$ in the above
example may be seen to be subdirectly irreducible, as its only ideals are
$\{0\}$, $Y$ and $E^{[0,a]}$.

\section{Algebraic Semantics}\label{sec:alg-sem}

In this section, we begin by rendering the Hilbert-style
systems of Section~\ref{sec:nine-logics} more tractable by
studying the derivability relation in $\ALm$ and its extensions.
We then give the semantics for the language $\lL$ in a pocrim
and show that the logics of Figure~\ref{fig:logics} are each sound
and complete for a corresponding class of pocrim. We use the
semantics to give a new proof that hoops form a variety.

\subsection{Derivability}

If $T$ is any subset of $\lL$, we say $B$ is {\em derivable} from $A$ in $T$
and write $A \ge_T B$, if $A \Lolly B$ is provable in $T$.  We say $A$ and $B$
are {\em equivalent} in $T$ and write $A \simeq_T B$, if $A \ge_T B$ and $B
\ge_T A$.  Thus $A \ge_{\ALm} B$ means that $A \Lolly B$ can be derived from
the axiom schemata $\axComp$, $\axComm$, $\axCurry$, $\axUncurry$ and $\axWk$ 
using {\it modus ponens}.
When the $T$ in question is clear from the context we
just write $\ge$ and $\simeq$.
Our goal is to find properties of these relations that make it easy
to prove facts such as $A \Lolly B \Lolly D \iAnd C \simeq_{T} B \iAnd A \Lolly C \iAnd D$,
where $T$ extends $\ALm$.

\begin{Lemma}\label{lma:ge-pre-order}
Let the theory $T$ extend $\ALm$. Then
$\ge_T$ is a pre-order and $\simeq_T$ is an equivalence relation.
\end{Lemma}
\Proof
Recall that a pre-order is a transitive and reflexive relation.
By definition, if $A \ge B$ and $B \ge C$, then $A \Lolly B$ and $B \Lolly C$
are both derivable in $\ALm$ and then using axiom {\axComp} and two applications of
{\it modus ponens}, we can derive $A \Lolly C$, so that $A \ge C$. So $\ge$
is transitive.
Now let $D$ be any provable formula, say the instance $V_1 \iAnd V_2 \Lolly V_1$
of {\axWk}. We can then derive $A \Lolly A$ as follows:
\begin{align*}
\mbox{1:}\quad& D \tag*{[\axWk]} \\
\mbox{2:}\quad& A \iAnd D \Lolly A \tag*{[\axWk]} \\
\mbox{3:}\quad& D \iAnd A \Lolly A \iAnd D \tag*{[\axComm]} \\
\mbox{4:}\quad& (A \iAnd D \Lolly A) \Lolly (D \iAnd A \Lolly A) \tag*{[3, \axComp]} \\
\mbox{5:}\quad& D \iAnd A \Lolly A  \tag*{[2, 4]} \\
\mbox{6:}\quad& D \Lolly A \Lolly A \tag*{[5, \axCurry]} \\
\mbox{7:}\quad& A \Lolly A & \tag*{[1, 6]} 
\end{align*}
(Here a justification such as [3, \axComp] indicates an application of modus
ponens with the result of line 3 as the cut-formula and an instance of {\axComp} as
the implication.)
So $\ge$ is reflexive and hence is indeed a pre-order.
That $\simeq$ is an equivalence relation follows immediately.
\Done

\begin{Lemma}\label{lma:ge-provability}
Let the theory $T$ extend $\ALm$.
For any formula $A$, the following are equivalent:
(i) $A$ is provable in $T$; (ii) $B \ge_T A$
for every formula $B$; (iii) $B \ge_T A$ for some
formula $B$ that is provable in $T$.
\end{Lemma}
\Proof
$\mbox{{\em(i)}} \Imp \mbox{{\em(ii)}}$:
If $A$ is provable and $B$ is any formula, then we can derive $B \Lolly A$
as follows. By assumption we have $A$. By $\axWk$ we have $A \iAnd B \Lolly A$, which by $\axCurry$ gives us $A \Lolly B \Lolly A$. Finally, from $A$ and the $A \Lolly B \Lolly A$ we obtain $B \Lolly A$. \\
$\mbox{{\em(ii)}} \Imp \mbox{{\em(iii)}}$:
This is trivial given that provable formulas exist.
\\
$\mbox{{\em(iii)}} \Imp \mbox{{\em(i)}}$:
By definition, if $B \ge A$, then $B \Lolly A$ is provable, so if
$B$ is provable, then $A$ follows with one application of {\it modus ponens}. \Done \\

In the sequel, as in the following
proof, we will often tacitly apply Lemma~\ref{lma:ge-provability}, typically
taking the provable formula $B$ in part {\em(iii)} to be $0 \equiv 1 \Lolly 1$ which
is provable by dint of Lemma~\ref{lma:ge-pre-order}.

\begin{Lemma}\label{lma:ge-congruence}
Let the theory $T$ extend $\ALm$.
With respect to the pre-order $\ge_{T}$,
$\Lolly$ is antimonotonic in its first argument and monotonic
in its second argument, while $\iAnd$ is monotonic in both arguments.
I.e., for any formulas $A$, $B$ and $C$ such that
$A \ge_{T} B$, the following hold:
\begin{align*}
B \Lolly C &\ge_{T} A \Lolly C \tag*{\em(i)} \\
C \Lolly A &\ge_{T} C \Lolly B  \tag*{\em(ii)} \\
A \iAnd C &\ge_{T} B \iAnd C \tag*{\em(iii)} \\
C \iAnd A &\ge_{T} C \iAnd B. \tag*{\em(iv)}
\end{align*}
The equivalence relation $\simeq_{T}$ is a congruence with respect to both $\Lolly$ and $\iAnd$.
I.e., for any formulas $A$, $B$ and $C$ such that $A \simeq_{T} B$, the following hold:
\begin{align*}
B \Lolly C &\simeq_{T} A \Lolly C \tag*{\em(v)} \\
C \Lolly A &\simeq_{T} C \Lolly B  \tag*{\em(vi)} \\
A \iAnd C &\simeq_{T} B \iAnd C \tag*{\em(vii)} \\
C \iAnd A &\simeq_{T} C \iAnd B. \tag*{\em(viii)}
\end{align*}
\end{Lemma}
\Proof
Assume that $A \ge B$, i.e, that $A \Lolly B$ is provable in $T$.
Using {\it modus ponens} and {\axComp}, we can derive $(B \Lolly C) \Lolly (A \Lolly C)$.
So {\em(i)} holds.
Using {\axComm}, {\axCurry}, {\axUncurry} and {\em(i)}, we have
$(X \Lolly Y \Lolly Z) \ge (Y \Lolly X \Lolly Z)$.
Instantiating $X$, $Y$ and $Z$ to $C \Lolly A$, $A \Lolly B$ and $C \Lolly B$
respectively, the left-hand side of this inequality becomes an instance
of {\axComp} and hence the right-hand side is provable in $T$. But the right-hand
side is exactly what we need to derive $(C \Lolly A) \Lolly (C \Lolly B)$
from our assumption $A \Lolly B$ using {\it modus ponens}. So {\em(ii)} holds.
We now have the following inequalities:
\begin{align*}
0 &\ge B \iAnd C \Lolly B \iAnd C \tag*{(Lemma~\ref{lma:ge-pre-order})}\\
  &\ge B \Lolly C \Lolly B \iAnd C
		\tag*{\axCurry}\\
  &\ge A \Lolly C \Lolly B \iAnd C
		\tag*{{\em(i)}}\\
  &\ge A \iAnd C \Lolly B \iAnd C \tag*{\axUncurry.}
\end{align*}
So {\em(iii)} holds and then {\em(iv)} follows using {\axComm}.
{\em(v)}, {\em(vi)}, {\em(vii)} and {\em(viii)} follow immediately
from the definition of $\simeq$, {\em(i)}, {\em(ii)}, {\em(iii)} and {\em(iv)}.
\Done

\begin{Lemma}\label{lma:tensor-assoc-0}
Let the theory $T$ extend $\ALm$.
For any formulas $A$, $B$ and $C$, the following hold:
\begin{align*}
(A \iAnd B) \iAnd C &\simeq_{T} A \iAnd (B \iAnd C) \\
A \iAnd 0 &\simeq_{T} A
\end{align*}
\end{Lemma}
\Proof
For any $D$, using $\axCurry$ and $\axUncurry$, we have:
\[
(A \iAnd B) \iAnd C \Lolly D \simeq A \Lolly B \Lolly C \Lolly D
\]
But we also have $B \Lolly C \Lolly D \simeq B \iAnd C \Lolly D$. Since $\simeq$ is a congruence, using $\axCurry$ and $\axUncurry$ again, we have
\[
A \Lolly B \Lolly C \Lolly D \simeq A \Lolly B \iAnd C \Lolly D \simeq
	A \iAnd (B \iAnd C) \Lolly D
\]
Taking $D$ to be $A \iAnd (B \iAnd C)$ and $(A \iAnd B) \iAnd C$,
we obtain $(A \iAnd B) \iAnd C \ge A \iAnd (B \iAnd C)$
and $A \iAnd (B \iAnd C) \ge (A \iAnd B) \iAnd C)$,
i.e., $(A \iAnd B) \iAnd C \simeq A \iAnd (B \iAnd C)$.
We leave the second part as an exercise.
%
%
\Done

\subsection{Semantics}\label{sec:semantics}

We now give a semantics for the language $\lL$ in which the semantic
values of formulas are elements of pocrims. It is convenient in describing
the semantics to work in a single language including the constant $1$.
The value of $1$ is only required to be an annihilator when that is
stated explicitly.

So given a pocrim $\VP$, and an assignment
$\alpha : \Var \cup \{1\} \To P$ of elements of $P$
to variables and the constant $1$,
we extend $\alpha$ to a meaning function $v_{\alpha} : {\lL} \To P$
by interpreting $\iAnd$ and $\Lolly$ as $+$ and
$\rImp$ respectively. So, for example, the formula $0$, i.e.,  $1 \Lolly 1$,
will be interpreted as $\alpha(1) \rImp \alpha(1)$, i.e.,  $0$, the identity
element of $\VP$.
We say that $\alpha$ {\em satisfies} a formula $A$,
if $v_{\alpha}(A) = 0$.
We say that $A$ is {\em valid} in $\VP$
if it is satisfied by every assignment $\alpha : \Var \cup \{1\} \To P$.
If $\VP$ is bounded with annihilator $1$, we say $A$ is {\em boundedly valid}
if it is satisfied by every assignment $\alpha : \Var \cup \{1\}\To P$
such that $\alpha(1) = 1$.
If $\cC$ is a class of pocrims that are not all bounded, we say a formula $A$ is {\em valid} 
in $\cC$,  if it is valid in every $\VP \in {\cC}$.
If $\cB$ is a class of bounded pocrims, we say a formula $A$ is {\em valid}
in $\cB$ if it is boundedly valid in every $\VP \in {\cB}$.
(This technical trick is convenient for the statement of the theorem that follows.)
A logic $L$ is {\em sound} for $\cC$ if every $A$ that is provable in $L$
is valid in $\cC$. $L$ is {\em complete} for $\cC$ if
every formula that is valid in $\cC$ is provable in $L$.
If $\VP$ is a pocrim, we write $\Th(\VP)$ for the set of all formulas that
are valid in $\VP$. In the proof of the following theorem, we exhibit
a pocrim $\VT$ such that $\Th(\VT)$ comprises precisely the set of
formulas that are provable in $\ALm$.

\begin{Theorem}\label{thm:sound-complete}
Each of our nine logics is sound and complete for the corresponding
class of pocrims shown in the following table:

\begin{center}
\begin{tabular}{|c|l|} \hline
$\ALm$ & pocrims \\ \hline
$\ALi$ & bounded pocrims \\ \hline
$\ALc$ & involutive pocrims \\ \hline
$\LLm$ & hoops \\ \hline
$\LLi$ & bounded hoops \\ \hline
$\LLc$ & involutive hoops \\ \hline
$\ML$ & idempotent hoops \\ \hline
$\IL$ & bounded idempotent hoops \\ \hline
$\BL$ & involutive idempotent hoops \\ \hline
\end{tabular}
\end{center}
\end{Theorem}
\Proof
We give the proof for $\ALm$. The modifications to give the proofs
for the other logics are straightforward.\\
For soundness, its suffices to show that all instances
of the axiom schemata used to define $\ALm$ are
valid and that {\it modus ponens} preserves validity.
We will just consider the axiom schema {\axCurry} and leave the rest
as an exercise.
For \axCurry, we have to show that
$v_{\alpha}((A \iAnd B \Lolly C) \Lolly (A \Lolly B \Lolly C)) = 0$
for any formulas $A$, $B$ and $C$ and any assignment $\alpha$ in any
pocrim.
Now
$v_{\alpha}((A \iAnd B \Lolly C) \Lolly (A \Lolly B \Lolly C)) = 
(v_{\alpha}(A) + v_{\alpha}(B) \rImp v_{\alpha}(C)) \rImp (v_{\alpha}(A) \rImp v_{\alpha}(B) \rImp v_{\alpha}(C))$.
Hence it is sufficient to show that every pocrim satisfies
$(a + b \rImp c) \rImp (a \rImp b \rImp c) = 0$,
i.e., that every pocrim satisfies $a + b \rImp c \ge a \rImp b \rImp c$.
Two applications of the residuation property (and some rearrangement using
the commutative monoid laws) show that this holds iff $(a + b) + (a + b \rImp c) \ge c$, which has the form $x + (x \rImp y) \ge y$.
By the residuation property, this is equivalent to $x \rImp y \ge x \rImp y$,
which holds since $\ge$ is a partial order, completing the proof that
all instances of $\axCurry$ are valid. \\
As for completeness,
by Lemmas~\ref{lma:ge-pre-order} and~\ref{lma:ge-congruence}, we
may define a structure $\VT = (T; 0, +, \rImp)$ by taking
$T = {\lL}/{\simeq}$ (the set of $\simeq$-equivalence classes) and
defining:
\begin{align*}
0 &= [0] \\
[A] + [B] &= [A \iAnd B] \\
[A] \rImp [B] &= [A \Lolly B].
\end{align*}
$\VT$ is the {\em term model} for $\ALm$.
It now follows using
\axComp, \axComm, \axCurry, {\axUncurry} and {\axWk}
and Lemmas~\ref{lma:ge-pre-order},~\ref{lma:ge-congruence}
and~\ref{lma:tensor-assoc-0} that $\VT$ is a pocrim.
Now as our axiom schemata are closed under substitution,
a formula $A$ is valid in $\VT$ iff it is valid under the
interpretation that maps each variable $P$ to $[P]$, i.e.,
iff $[A] = 0$, which holds iff $A$ is provable in $\ALm$.
Completeness follows, since if $A$ is valid in all pocrims, then
it is certainly valid in $\VT$.
\Done

As we have defined it, the class of pocrims is a quasivariety, i.e.,
its defining properties are Horn clauses over equational atoms.
It is known that this is the best we can do: the class
of involutive pocrims cannot be characterised by equational laws.
Since involutive pocrims are
characterised over bounded pocrims and over pocrims by equational laws, it
follows that the class of pocrims and the class of bounded pocrims are also not
varieties.  See~\cite{raftery07} and the works cited therein for these results.

In \cite{arthan-oliva14a} we present a number of proofs derived from
machine-oriented derivations found by the automated theorem-prover Prover9.  Our use of
Prover9 relies heavily on the fact that the class of hoops is actually a
variety with quite a short and simple equational axiomatisation.  In a {\it
tour de force} of equational reasoning, Bosbach~\cite{bosbach69a} gave a direct
proof of an equational axiomatization of the class of hoops. Using
Theorem~\ref{thm:sound-complete}, we can give a new proof that hoops
form a variety by showing
how to transform a proof of a formula $A$ in $\LLm$
into an equational proof
that $a = 0$, where $a$ is a translation into the language of pocrims
of the formula $A$.

\begin{Theorem}\label{thm:luk-provability-equational}
Let $\VH$ be a structure for the signature $(0, +, \rImp)$.
\begin{description}
\item[I.]
$\VH$ is a hoop
iff $(H; 0, +)$ is a commutative monoid and $\VH$ satisfies the following equations:

\newcounter{temp}
\begin{enumerate}
\item\label{luk-imp-self} $x \rImp x = 0$
\item\label{luk-imp-zero} $x \rImp 0 = 0$
\item\label{luk-conj-imp} $x + y \rImp z = x \rImp y \rImp z$
\item\label{luk-cwc} $x + (x \rImp y) = y + (y \rImp x)$
\setcounter{temp}{\value{enumi}}
\end{enumerate}

\item[II.]
$\VH$ is a bounded hoop iff it satisfies the above equations and also:
\begin{enumerate}
\setcounter{enumi}{\value{temp}}
\item\label{luk-one-imp} $1 \rImp x = 0$.
\end{enumerate}
\end{description}
\end{Theorem}
\Proof
{\bf I}:
It follows easily from the definitions (or from
Theorem~\ref{thm:sound-complete}) that equations~\ref{luk-imp-self}
to~\ref{luk-cwc} hold in any hoop.
For the converse, Theorem~\ref{thm:sound-complete}
implies that it is sufficient to show that, if there is
proof of $A$ in $\LLm$
then $[A]$, (the element of the term model of $\LLm$
represented by $A$) can be reduced to 0 using the commutative
monoid laws and equations~\ref{luk-imp-self} to \ref{luk-cwc}.
It follows that if $a$ is the formula obtained from $A$ by
replacing $\iAnd$ and $\Lolly$ by $+$ and $\rImp$ respectively, then
$A$ is provable in $\LLm$ iff every hoop satisfies $a = 0$.
We will show how to translate a proof of $A$ in $\LLm$
into a sequence of equations $a = a_1 = \ldots = a_n = 0$,
where each equation $a_i = a_{i+1}$ is obtained
by applying one of the equations~\ref{luk-imp-self} to \ref{luk-cwc}
to a subterm of $a_i$ or $a_{i+1}$ or by using the commutative monoid laws.
The equational derivation is defined by recursion over a proof constructed
using {\it modus ponens} from the axiom schemata
$\axComp$, $\axComm$, $\axCurry$, $\axUncurry$ and $\axWk$.
We give the justification as we define each step in the derivation,
so once the definition is complete, the proof is complete.

\noindent
$\axComp$: we want $(a \rImp b) \rImp (b \rImp c) \rImp (a \rImp c) = 0$ for arbitrary $a$, $b$ and $c$:
\begin{align*}
(a \rImp b) \rImp (b \rImp c) \rImp (a \rImp c) &=
   (a \rImp b) + (b \rImp c) + a \rImp c &
      & \tag*{$2 \times \mbox{(eq. \ref{luk-conj-imp})}$} \\
   &= a + (a \rImp b) + (b \rImp c) \rImp c
     & \tag*{(comm. monoid)}\\
   &= b + (b \rImp a) + (b \rImp c) \rImp c
     & \tag*{(eq. \ref{luk-cwc})} \\
   &= (b \rImp a) + b + (b \rImp c) \rImp c
     & \tag*{(comm. monoid)}\\
   &= (b \rImp a) + c + (c \rImp b) \rImp c
     & \tag*{(eq. \ref{luk-cwc})} \\
   &= (b \rImp a) + (c \rImp b) + c \rImp c
     & \tag*{(comm. monoid)} \\
   &= (b \rImp a) + (c \rImp b) \rImp c \rImp c
     & \tag*{(eq. \ref{luk-conj-imp})} \\
   &= (b \rImp a) + (c \rImp b) \rImp 0
     & \tag*{(eq. \ref{luk-imp-self})} \\
   &= 0.
     & \tag*{(eq. \ref{luk-imp-zero})} 
\end{align*}

\noindent
$\axComm$: we want $a + b \rImp b + a = 0$ for arbitrary $a$ and $b$:
\begin{align*}
a + b \rImp b + a &= a + b \rImp a + b
     & \tag*{(comm. monoid)} \\
   &= 0 & \tag*{(eq. \ref{luk-imp-self})} 
\end{align*}

\noindent
$\axCurry$: we want $(a + b \rImp c) \rImp (a \rImp b \rImp c) = 0$ for arbitrary $a$, $b$ and $c$:
\begin{align*}
(a + b \rImp c) \rImp (a \rImp b \rImp c) &=
   (a \rImp b \rImp c) \rImp (a \rImp b \rImp c)
     & \tag*{(eq. \ref{luk-conj-imp})} \\
   &= 0.
     & \tag*{(eq. \ref{luk-imp-self})} 
\end{align*}

\noindent
$\axUncurry$: we want $(a \rImp b \rImp c) \rImp (a + b \rImp c) = 0$ for arbitrary $a$, $b$ and $c$:
\begin{align*}
(a \rImp b \rImp c) \rImp (a + b \rImp c) &=
   (a \rImp b \rImp c) \rImp (a \rImp b \rImp c)
     & \tag*{(eq. \ref{luk-conj-imp})} \\
   &= 0.
     & \tag*{(eq. \ref{luk-imp-self})} 
\end{align*}

\noindent
$\axWk$: we want $a + b \rImp a = 0$ for arbitrary $a$ and $b$:
\begin{align*}
a + b \rImp a &=
   b + a \rImp a
     & \tag*{(comm. monoid)} \\
   &= b \rImp a \rImp a
     & \tag*{(eq. \ref{luk-conj-imp})} \\
   &= b \rImp 0
     & \tag*{(eq. \ref{luk-imp-self})} \\
   &= 0.
     & \tag*{(eq. \ref{luk-imp-zero})} 
\end{align*}

\noindent
{\it Modus ponens}: we are given $a = 0$ and $a \rImp b = 0$ and we want $b = 0$:
\begin{align*}
b &=
   b + 0
     & \tag*{(comm. monoid)} \\
   &= b + (b \rImp 0)
     & \tag*{(eq. \ref{luk-imp-zero})} \\
   &= 0 + (0 \rImp b)
     & \tag*{(eq. \ref{luk-cwc})} \\
   &= 0 \rImp b
     & \tag*{(comm. monoid)} \\
   &= a \rImp b
     & \tag*{(given)} \\
   &= 0.
     & \tag*{(given)} 
\end{align*}
\noindent This completes the recursive definition. \\
{\bf II:} Like part {\bf I}, using equation
(\ref{luk-one-imp}) to translate the axiom schema \axEFQ. \Done \\

We have chosen equations~\ref{luk-imp-self} to~\ref{luk-one-imp} for
convenience in the above proof.  In fact, our Prover9 work uses $x + 1 = 1$ to
characterize $1$ rather than  $1 \rImp x = 0$.
The reader may enjoy showing that the two are equivalent.
A more intricate exercise is to show that equation~\ref{luk-imp-zero}
is redundant, as it can be derived from the commutative monoid laws and
equations~\ref{luk-imp-self},~\ref{luk-conj-imp} and~\ref{luk-cwc}.

%
%
%
%
%
%

In \cite{arthan-oliva14a} we give a sequent calculus
that is equivalent to our Hilbert-style presentation of $\LLm$.
The sequent calculus proofs can also be translated into
equational proofs along similar lines to the translation
given above.

\section{Identities in Hoops}\label{sec:indirect-method}

Blok and Ferreirim proved that the quasi-equational
theory of hoops is decidable
\cite{blok-ferreirim00}. However, the proof does not lead to any bounds
on the complexity of the decision procedure. More recently,
Bova and Montagna \cite{Bova:2009} have shown that the quasi-equational
theory of commutative GBL-algebras is in PSPACE and, in fact, is
PSPACE-complete. It can be shown that the quasi-equational theory
of commutative GBL-algebras is a conservative extension of that of
hoops and hence Bova and Montagna's work implies that the quasi-equational
theory of hoops is in PSPACE.
Bova and Montagna use a generalisation of the ordinal sum construction.
This {\em poset sum} construction takes as input a family of commutative
GBL-algebras $\VG_p$ indexed by a poset $\VP$.
They show that a quasi-equation involving $n$ symbol occurrences
holds in all commutative GBL-algebras iff it holds in all
finite algebras of size at most $2^{3n^2}$ that are the poset sum
of a family of finite MV-chains\footnote{
Finite MV-chains are the MV-algebras corresponding to the hoops
$\VL_i$ of Example~\ref{eg:vl}.
}
indexed by a poset comprising a tree of height at most $n$ and with at most
$2^{n^2}$ nodes.
They then give an ingenious non-deterministic
algorithm that checks using polynomial space whether a given
quasi-equation can be refuted in the corresponding set of finite
algebras. Since co-NPSPACE, NPSPACE and PSPACE coincide, this shows
that the quasi-equational theory is in PSPACE.

Unfortunately, the algorithm of Bova and Montagna is infeasible,
certainly for hand calculation even on small examples: the valid
equation $\Not(\Not\Not x \rImp x) = 0$ contains 8 symbols and
just the number of trees to be considered would be enormous.
As we are interested in verifying certain specific
identities, we need a more practical method.
To this end, we will show that a identity is valid in all hoops iff
it is valid in a restricted class of hoops enjoying some very convenient
algebraic properties.  This does not provide a decision procedure, but it does
provide a quick indirect method of proof for many important identities.
We begin with a review of the decision problem for identities in
involutive hoops, which we often need to consider when applying the
indirect method for bounded hoops in general.

\subsection{Identities in involutive hoops}\label{sec:id-inv-hoops}

By contrast with the case for general hoops, the equational theory of
involutive hoops can decided by a computationally efficient reduction to
(linear) real arithmetic. In this section, we review the proof of this
result, which relies on the fact that involutive hoops are
definitionally equivalent to MV-algebras.  This definitional equivalence
is stated without proof in \cite{blok-ferreirim00}. We give the proof
here for expository purposes and because it involves some identities
that will be useful later.

MV-algebras were originally
introduced by Chang~\cite{Chang58b} and have been widely studied.  We
adopt the definition and notation of~\cite{Cignoli-et-al00}:

\begin{Definition}\label{def:mv-algebras}
An MV-algebra is a structure for the
signature $(0, \oplus, \Not)$ whose $(0, \oplus)$-reduct is a commutative
monoid and which, with $x \ominus y$ defined as $\Not(\Not x \oplus y)$, satisfies the following identities:
\begin{align*}
   x \oplus \Not 0 &= \Not 0 \\
   \Not\Not x &= x \\
   (y \ominus x) \oplus x &= (x \ominus y) \oplus y
\end{align*}
\end{Definition}

\begin{Lemma}\label{lma:inv-hoop-identities}
If $\VH$ is an involutive hoop, then $\VH$ satisfies:
\[
x \rImp y = \Not(x + \Not y) \quad \quad \quad
x + y = \Not(x \rImp \Not y).
\]
\end{Lemma}
\Proof
Recalling that $\Not x = x \rImp 1$ by definition
and using {\INV} we have:
\[
x \rImp y = x \rImp \Not\Not y = x \rImp \Not y \rImp 1 = x + \Not y \rImp 1 = \Not(x + \Not y),
\]
and then we have
$
\Not(x \rImp \Not y) = \Not\Not(x + \Not\Not y) = x + y.
$
\Done

\begin{Theorem}\label{thm:hoops-equiv-mv}
The variety of involutive hoops and the variety of MV-algebras are
definitionally equivalent.
\end{Theorem}
\Proof
Let $\VH$ be an involutive hoop and define $x \oplus y = x + y$ and $\Not x = x \rImp 1$.
Then $(0, \oplus)$
is a commutative monoid and we have $x \oplus \Not 0 = x + 1  = 1 = \Not 0$ and $\Not\Not x = x$.
Moreover, by
Lemma~\ref{lma:inv-hoop-identities}, $x \ominus y = \Not(\Not x \oplus y) =  y \rImp x$ and hence $(y \ominus x) \oplus x = x + (x \rImp y) = y + (y \rImp x) = (x \ominus y) \oplus y$.
Thus $(H; 0, \oplus, \Not)$ is an MV-algebra. 
Conversely, let $\VM$ be an MV-algebra and define $1 = \Not 0$,
$x + y = x \oplus y$
and $x \rImp y = \Not(x + \Not y)$.
Then certainly $(M; 0, +)$ is a commutative monoid.
We have $x \rImp 1 = \Not(x + \Not\Not 0) = \Not x$.
Hence $x + \Not x = x + (x \rImp 1) = 1 + (1 \rImp x) = 1$,
so that $x \rImp x = \Not(x + \Not x) = \Not 1 = 0$.
Hence equation~\ref{luk-imp-self} in the equational characterization
of bounded hoops given in
Theorem~\ref{thm:luk-provability-equational} is satisfied.
The other equations in that characterization are easily verified
and so, as we have $\Not\Not x = x$,
$(M; 0, +, \rImp)$ is an involutive hoop.
\Done \\

Chang~\cite{Chang58b,Chang59} showed that the system
$\Luk$ of Section~\ref{sec:nine-logics}
is sound and complete for the class of MV-algebras under a
semantics which corresponds to our semantics
for hoops under the equivalence of Theorem~\ref{thm:hoops-equiv-mv}.
Since $\LLc$ is sound and complete for involutive hoops,
it follows that $\LLc$ and $\Luk$ are equivalent.
Chang's work also implies that an identity holds in all MV-algebras iff it holds in
the MV-algebra corresponding to the involutive hoop $\UI$.
(See~\cite{Cignoli-et-al00} for more information on MV-algebras.)

\begin{Theorem}\label{thm:llc-decidable}
An identity $s = t$ holds in all involutive hoops iff
it holds in the hoop $\UI$ of Example~\ref{eg:ui}.
Hence the equational theory of involutive hoops is decidable.
\end{Theorem}
\Proof The first claim follows from the remarks above about the
equivalence between involutive hoops and MV-algebras and the definition of $\UI$.
Given the first claim, to decide $s = t$, use the formulas
for the operations on $\UI$ given in Example~\ref{eg:ui}
to translate $s = t$ into a formula in the
language of real arithmetic, treating $\Max$ and $\Min$ as
abbreviations:
$\phi(\Max(x, y)) \equiv (x \ge y \And \phi(x)) \Or (x < y \And \phi(y))$ and
$\phi(\Min(x, y)) \equiv (x \ge y \And \phi(y)) \Or (x < y \And \phi(x))$.
Equality of the translated formula may then be decided by the well-known decision procedures for (linear) real arithmetic.
\Done


It can be shown, using results of Blok and Ferreirim
\cite{blok-ferreirim00}, that an identity holds in all Wajsberg hoops
iff it holds in the hoop $\RRnn$ of Example~\ref{eg:rrnn}. Hence the
equational theory of Wajsberg hoops also reduces to (linear) real
arithmetic.
\subsection{Identities in general hoops}\label{sec:characterizing-identities}

We now give our indirect method for proving identities in general hoops.
The method is based on the characterization of subdirectly irreducible
hoops due to Blok and Ferreirim \cite{blok-ferreirim00}.  They proved
that a hoop $\VH$ is subdirectly irreducible iff it is isomorphic to an
ordinal sum $\VS \ordSum \VF$ where $\VS$ is subdirectly irreducible,
totally ordered  and Wajsberg and where $\VS$ is trivial iff $\VH$ is
trivial.  $\VS$ and $\VF$ are uniquely determined by these conditions
and are called the {\em support} and the {\em fixed subhoop} of $\VH$
respectively.

The following theorem is really two: one for bounded hoops and one for all hoops.
From now on, when we work with bounded hoops, we will take the annihilator $1$ as part of the signature, so that homomorphisms must preserve it and it must
be included when we consider the bounded subhoop of a given bounded hoop $\VH$
generated by some subset of $H$.

\begin{Theorem}\label{thm:characterizing-identities}
Let $\phi$ be an identity in the language of a (bounded) hoop in
the variables $x_1, \ldots, x_n$.
Then $\all{x_1, \ldots, x_n} \phi$ is valid in the class of all (bounded) hoops
iff $\phi(x_1, \ldots, x_n)$ holds
under any interpretation of $x_1, \ldots, x_n$ in a (bounded) hoop $\VH$
that can be expressed as an ordinal sum
$\VS \ordSum \VF$ where $\VS$ is subdirectly irreducible and Wajsberg,
where $\VH$ is generated by $x_1, \ldots, x_n$ and
where $S = \{0\}$ iff $H = \{0\}$.
\end{Theorem}
\Proof
$\Rightarrow$:
if $\all{x_1, \ldots, x_n} \phi$ is valid in the class of all (bounded) hoops,
then $\phi(x_1, \ldots, x_n)$ holds for any interpretation of $x_1, \ldots, x_n$
in any (bounded) hoop.  \\
$\Leftarrow$:
Assume that  $\phi(x_1, \ldots, x_n)$ holds
in any (bounded) hoop $\VH$
satisfying the stated conditions on $\VH$ and on the interpretation of $x_1, \ldots, x_n$.
Let then $\VH$ be an arbitrary (bounded) hoop and let
$x_1, \ldots, x_n \in H$. We must prove that
$\phi(x_1, \ldots, x_n)$ holds in $\VH$.
Clearly $\phi(x_1, \ldots, x_n)$ holds
in $\VH$ iff it holds in the subhoop of $\VH$ generated by $x_1, \ldots, x_n$.
So we may assume $\VH$ is generated by $x_1, \ldots, x_n$.
By a classic result of Birkhoff (e.g., see
\cite[Theorem II.8.6]{burris-sankappanavar81}) $\VH$ is isomorphic to a subdirect
product of subdirectly irreducible hoops each of which is
a homomorphic image of $\VH$ (and hence is generated by the images
$[x_1], \ldots, [x_n]$ of $x_1, \ldots, x_n$).
The identity $\phi(x_1, \ldots, x_n)$ holds in the
subdirect product if $\phi([x_1], \ldots, [x_n])$ holds in each factor.
So we may assume that $\VH$ is subdirectly irreducible and generated by $x_1, \ldots, x_n$.
By the theorem of Blok and Ferreirim, $\VH$ has subhoops $\VS$ and $\VF$
such that $\VS$ is subdirectly irreducible and Wajsberg and
$\VH \cong \VS \ordSum \VF$. Moreover $S = \{0\}$ iff $H = \{0\}$.
Hence, by assumption, $\phi(x_1, \ldots, x_n)$ holds in $\VH$.
\Done \\


By the definition of the ordinal sum, if $s \in S \Diff \{0\}$ and $f \in F
\Diff \{0\}$, then $(f \rImp s) \rImp s = 0 \rImp s = s \not= (s \rImp f) \rImp
f = f \rImp f = 0$, i.e., $s$ and $f$ do not satisfy the Wajsberg condition.
So, if $\VH$ is subdirectly irreducible and Wajsberg, then it is equal to its
support $\VS$ and hence is totally ordered.

\begin{figure}[t]
\begin{center}
\begin{tabular}{|p{0.9\hsize}|}
\hline
To verify $\phi(x_1, \ldots, x_n)$ in all hoops, verify it in the following cases:\\
{\em Case {\em(i)}: $\VH \cong \VS \ordSum \VF$ with $F = \{ 0 \}$. $\VH \cong \VS$.} \\
{\em Case {\em(ii)}: $\VH \cong \VS \ordSum \VF$ with $F \neq \{ 0 \}$. There is a  subcase for each choice of $I = \{i \ST x_i \in S\} \neq \emptyset$ and $J = \{j \ST x_j \in F\} \neq \emptyset$, with $\VF$ generated by the $x_j$ with $j \in J$.}\\
{\em In both cases $\VS$ is subdirectly irreducible, Wajsberg and generated by the $x_i \in S$.}  \\
\hline
\end{tabular}
\end{center}
\vspace{-3mm}
\caption{Template for applying theorem \ref{thm:characterizing-identities} to all hoops} 
\label{fig:unbounded-case}
\end{figure}

The rest of Section~\ref{sec:indirect-method} illustrates the use of
Theorem \ref{thm:characterizing-identities}.  To prove
that an identity $\phi(x_1, \ldots, x_n)$ holds in all hoops, we follow the
template of Figure \ref{fig:unbounded-case}.  We consider an interpretation of
the variables $x_1, \ldots, x_n$ in a hoop $\VH = \VS \ordSum \VF$ satisfying
the stated conditions. Apart from the trivial case when $n = 0$, the set $I =
\{i \ST x_i \in S\}$ cannot be empty (otherwise we would have $S = \{0\}$ while
$H = F \not= \{0\}$). We then consider all possible cases for the set $I$.  So
let $J = \{1, \ldots, n\} \Diff I$.  If $J = \emptyset$, then we must verify
that $\phi$ holds in a subdirectly irreducible Wajsberg hoop $\VS$ (Case {\em
(i)}).  If $J \not= \emptyset$, we have to verify that $\phi$ holds in $\VS
\ordSum \VF$, where $\VS$ is generated by the $x_i$ with $i \in I$, $\VF$ is
generated by the $x_j$ with $j \in J$ and $x_j \not= 0$ for $j \in J$ (Case
{\em (ii)}). The cases where $J \not= \emptyset$  are often easy to verify
using identities such as $x_j \rImp x_i = 0$ and $x_i \rImp x_j = x_j$ when $i
\in I$ and $j \in J$ that follow from the definition of the ordinal sum.

To prove an identity holds
in all bounded hoops, we follow the template of Figure \ref{fig:bounded-case}.
We have the same cases as for the unbounded case with the extra assumption that
$\VH$ is bounded, and hence involutive in Case {\em (i)}.  We must also
consider the possibility that $J = \emptyset$ and $F$ is generated by the
constant $1$, i.e.  $\VH \cong \VS \ordSum \BB$ where $\VS$ is subdirectly
irreducible and Wajsberg and $\BB$ is the boolean hoop with $B = \{0, 1\}$
(Case {\em (iii)}).

\begin{figure}[h]
\begin{center}
\begin{tabular}{|p{0.9\hsize}|}
\hline
To verify $\phi(x_1, \ldots, x_n)$ in all bounded hoops, verify it in the following cases:\\
{\em Cases {\em(i)} and {\em(ii)}:
As in Figure \ref{fig:unbounded-case}, with the extra assumption that $\VH$ is bounded.} \\
{\em Case {\em (iii)}: $\VH \cong \VS \ordSum \BB$, with all $x_i \in S$.} \\
{\em In all cases, $\VS$ is subdirectly irreducible, Wajsberg and generated by the $x_i \in S$.} \\
\hline
\end{tabular}
\end{center}
\vspace{-3mm}
\caption{Template for applying theorem \ref{thm:characterizing-identities} to bounded hoops} 
\label{fig:bounded-case}
\end{figure}


The structure of free MV-algebras and hence of free involutive hoops is quite
well understood. See~\cite{Cignoli-et-al00} for a good account of this topic.
Very little is known about free hoops or free bounded hoops. It can be shown
that the free bounded hoop on one generator is a subdirect product of hoops
isomorphic to subhoops of $\UI$, $\UI \ordSum \BB$ and $\RR^{{\ge}0} \ordSum
\BB$. This suggested the identity of the following example.

\begin{Example}
If $0 < k \in \NN$, the identity $\Not kx \rImp \Not\Not x \rImp x = 0$ clearly
holds in any involutive hoop. It also holds in any hoop of the form $\VS
\ordSum \BB$ (since in such a hoop, either $x = 1$ or $\Not kx = 1$).  This
covers cases~{\em(i)} and~{\em(iii)} in the template of
Figure~\ref{fig:bounded-case}.  As the identity has only one variable, there is
nothing to prove in case~{\em(ii)}.  Hence, $\Not kx \rImp \Not\Not x \rImp x =
0$ holds in any bounded hoop. In~\cite{arthan-oliva14a}, we demonstrate how to
construct a proof in $\LLi$ of the formula of $\lL$ that corresponds to this
identity.  Unfolding the constructions used in the inductive step of this
demonstration reveals 19 intricate applications of the axiom $\axCWC$.
\end{Example}

\subsection{Application: de Morgan identities}

In any bounded pocrim, we have $\Not(x + y) = x + y \rImp 1 = x \rImp y \rImp 1 = x \rImp \Not y$, so the first identity in the following theorem is easily
proved.
In a bounded hoop, we have a kind of dual identity: $\Not(x \rImp y) = \Not\Not x + \Not y$.
As can be seen in~\cite{arthan-oliva14a}, the simplest known elementary proof of
the dual identity is quite involved.
The indirect proof using Theorem~\ref{thm:characterizing-identities} is much shorter:
\begin{Theorem}\label{thm:de-morgan}
The following identities are satisfied in any bounded hoop:
\[
\Not(x + y) = x \rImp \Not y \quad \quad \quad
\Not(x \rImp y) = \Not\Not x + \Not y
\]
\end{Theorem}
\Proof See the above remarks for the first identity. For the second we follow
the template of Figure \ref{fig:bounded-case}, requiring us to prove
the identity in the following cases for a hoop $\VH$ and its elements $x$ and $y$:
\\
Case {\em(i)}: Our assumptions imply that $\VH$ is involutive, hence
by Lemma~\ref{lma:inv-hoop-identities}
$$
\Not(x \rImp y) = \Not\Not(x + \Not y) = x + \Not y = \Not\Not x + \Not y
$$
Case {\em(ii)}:
$\VH = \VS \ordSum \VF$,
$\{x, y\}\cap S \not= \emptyset$,
$\{x, y\}\cap F \Diff \{0\} \not= \emptyset$:
this leads to two subcases that are proved 
using elementary properties of $\VS \ordSum \VF$, as follows:
\\
Subcase {\em(ii)(a)}:  $x \in S$, $y \in F \Diff \{0\}$:
$$
\Not(x \rImp y) = \Not y = 0 + \Not y = \Not\Not x + \Not y
$$
Subcase {\em(ii)(b)}: $x \in F \Diff \{0\}$, $y \in S$:
\begin{align*}
\Not(x \rImp y) = \Not 0 = 1 = \Not\Not x + 1 = \Not\Not x + \Not y.
\end{align*}
Case {\em(iii)}: $\VH = \VS \ordSum \BB$ where $x, y \in S$: 
for $u \in S$, $\Not u = 1$, so as $x \rImp y \in S$, we have $\Not(x \rImp y) = 1 = 0 + 1 = \Not\Not x + \Not y$. \Done

\subsection{Application: The Ferreirim-Veroff-Spinks theorem}

Ferreirim \cite{Ferreirim92} proved by indirect methods that if $e$ is an
idempotent element in a $k$-potent hoop, then the mapping $x \mapsto e \rImp x$
is an additive homomorphism.  Using Otter \cite{McCune03}, Veroff and Spinks
\cite{Veroff-Spinks04} found a syntactic proof without assuming $k$-potency.  A
simplified and more abstract presentation of their proof is given by the
present authors in~\cite{arthan-oliva14a}.  When we express the theorem as an
identity and apply Theorem~\ref{thm:characterizing-identities}, it turns out
that it is only the case when the hoop is subdirectly irreducible and Wajsberg
that presents any difficulties:

\begin{Theorem}
The following identity holds in any hoop:
$$
(e \rImp e + e) \rImp (e \rImp x + y) \rImp (e \rImp x) + (e \rImp y) = 0.
$$
\end{Theorem}
\Proof
Writing $a = e \rImp e + e$, $b = e \rImp x + y$ and $c = (e \rImp x) + (e
\rImp y)$, what we have to prove is that $a \rImp b \rImp c = 0$.
According to the template of Figure~\ref{fig:unbounded-case},
it is sufficient to verify $a \rImp b \rImp c = 0$ in the following cases for a
hoop $\VH$ and its elements $e$, $x$ and $y$:
\\
Case {\em(i)}: {\em $\VH$ subdirectly irreducible and Wajsberg}:
In this case, we claim that $a + b \ge c$, whence $a \rImp b \rImp c = 0$
as required.
By Lemma~\ref{lma:idemp-wajsberg}, there are two subcases:
\\
Subcase {\em(i)(a)}: $e \rImp e + e = e$: we have $a = e$, so that:
$$
a + b =  e + (e \rImp x + y) \ge x + y \ge (e \rImp x) + (e \rImp y) = c.
$$
Subcase {\em(i)(b)}: {\em $\VH$ is bounded and $e + e = 1$}: we have $a = e \rImp 1$, so that:
$$
a + b =  (e \rImp 1) + (e \rImp x + y) \ge (e \rImp x) + (e \rImp y) = c.
$$
Case {\em(ii)}:
$\VH = \VS \ordSum \VF$,
$\{e, x, y\}\cap S \not= \emptyset$,
$\{e, x, y\}\cap F \Diff \{0\} \not= \emptyset$:
since the equation is symmetric in $x$ and $y$, this leads to 4 subcases.
In each of these subcases, we claim that $b = c$, whence $a \rImp b \rImp c = 0$
as required. The claim is verified using elementary properties of $\VS \ordSum \VF$ as follows:
\\
Subcase {\em(ii)(a)}: $e \in S, x, y \in F \Diff \{0\}$:
$$
b = e \rImp x + y = x + y = (e \rImp x) + (e \rImp y) = c.
$$
Subcase {\em(ii)(b)}: $e, x \in S, y \in F \Diff \{0\}$:
$$
b = e \rImp x + y = e \rImp y = y =  (e \rImp x) + y
   = (e \rImp x) + (e \rImp y) = c.
$$
Subcase {\em(ii)(c)}: $e, x \in F \Diff \{0\}, y \in S$:
$$
b = e \rImp x + y = e \rImp x = (e \rImp x) + 0 = (e \rImp x) + (e \rImp y) = c.
$$
Subcase {\em(ii)(d)}: $e \in F \Diff \{0\}, x, y \in S$:
\begin{align*}
b = e \rImp x + y = 0 = 0 + 0 = (e \rImp x) + (e \rImp y) = c. \tag*{\Done}
\end{align*}

\begin{Lemma}\label{lma:idemp-wajsberg}
If $e$ is any element of a
totally ordered  Wajsberg hoop $\VH$ then either $e \rImp e + e = e$
or $\VH$ is bounded and $e + e = 1$.
\end{Lemma}
\Proof
We have $e + e = e + e + (e + e \rImp e) = e + (e \rImp e + e)$, so the
lemma follows from the fact that totally ordered irreducible Wajsberg hoops are
{\em semi-cancellative}, i.e., if $a + b = a + c$ with $b \not= c$, then
$\VH$ is bounded and $a + b = 1$.
This can be extracted from the characterization of subdirectly
irreducible hoops: see \cite[Theorem 26, part
{\em(ix)}]{arthan-oliva12}.  (The theorem we state there is for a special class
of hoops called coops, but the proof of that part of the theorem goes through
in the same way for a general hoop.)
However, Ferreirim \cite[Lemma 4.5]{Ferreirim92} gives the following
neat elementary proof: assume $a + b = a + c$ is not an annihilator, so there is
$u \in H$ with $u > a + b$. Since $u > a + b$,
$b > a \rImp u$ is impossible, so, as $\VH$ is totally ordered,
we must have $a \rImp u \ge b$, i.e., $(a \rImp u) \rImp b = 0$.
But then, as $\VH$ is Wajsberg, we have:
\begin{align*}
b &= ((a \rImp u) \rImp b) \rImp b \\
  &= (b \rImp a \rImp u) \rImp a \rImp u \\
  &= a \rImp (a + b \rImp u) \rImp u \\
  &= a \rImp (u \rImp a + b) \rImp a + b \\
  &= a \rImp a + b
\end{align*}
As $a + b = a + c$, the same argument gives us that
$c = a \rImp a + c$, so $b = c$.
\Done

\subsection{Application: the idempotent subhoop}
The set of idempotent elements of a pocrim is clearly closed under $+$ but it
need not be closed under $\rImp$. For example, in the pocrim $\VQ_4$ of
Example~\ref{eg:q-four}, $u$ and $1$ are idempotent, but $u \rImp 1 = v \not= 1$ and $v
+ v = 1$.  In a hoop, however, Theorem~\ref{thm:characterizing-identities}
enables us to prove the following somewhat surprising theorem:
\begin{Theorem}
Let $\VH$ be a hoop. The set $J = \{ x : H \ST x + x = x \} $
of idempotent elements of $\VH$ is the universe of a subhoop.
\end{Theorem}
\Proof
We must show that $0 \in J$, $J + J \subseteq J$ and $J \rImp J \subseteq J$.
The first two assertions are easy. As for $J \rImp J \subseteq J$,
let us define $i : H \To H$ by $i(x) = x \rImp x + x$,
so that $x \in J$ iff $i(x) = 0$.
It is sufficient to show that the following identity holds:
\begin{align*}
i(x) \rImp i(y) \rImp i(x \rImp y) &= 0 \tag*{($*$)}
\end{align*}
According to the template of Figure~\ref{fig:unbounded-case},
it is sufficient to verify ($*$) in the following cases for a
hoop $\VH$ and its elements $x$ and $y$:
\\
Case {\em(i)}: {\em $\VH$ subdirectly irreducible and Wajsberg}:
By Lemma~\ref{lma:idemp-wajsberg}, if $\VH$ is not bounded,
then $i(x) = x$ for all $x$ and ($*$) is trivial.
If $\VH$ is bounded, then it is involutive and by Theorem~\ref{thm:llc-decidable},
it is enough to verify ($*$) in $\UI$.
Now in $\UI$, we have:
\[
i(x) = \left\{
   \begin{array}{l@{\quad}l}
      x     & \mbox{if $x \le \frac{1}{2}$} \\[2mm]
      1 - x & \mbox{if $x \ge \frac{1}{2}$.}
   \end{array}
\right.
\]
($*$) holds in any hoop if $x \ge y$, so, in $\UI$, we may assume $x < y$,
so that $x \rImp y = y - x$. Then ($*$) holds iff
$i(x) + i(y) \ge i(y - x)$ and we have eight cases
for $x, y, y - x \in [0, 1]$ as follows:
{\def\Y{\checkmark}
\def\N{\times}
\def\X{\maltese}
\[
\begin{array}{cccr@{}c@{}l}
x > \frac{1}{2} & y > \frac{1}{2} & y - x > \frac{1}{2}
          &   i(x) + i(y) &{\relax}\ge{\relax}& i(y - x) \\
               \hline
\Y & \Y & \Y &                &\X               \\
\Y & \Y & \N &  1 - x + 1 - y &\ge& y - x       \\
\Y & \N & \Y &                &\X               \\
\Y & \N & \N &                &\X               \\
\N & \Y & \Y &      x + 1 - y &\ge& 1 - (y - x) \\
\N & \Y & \N &      x + 1 - y &\ge& y - x       \\
\N & \N & \Y &                &\X               \\ 
\N & \N & \N &          x + y &\ge& y - x
\end{array}
\]
The cases marked $\X$ are impossible, as the constraints on
$x$, $y$ and $y - x$ are inconsistent.
In the other cases, the inequalities are easily verified using the constraints.
That completes case {\em(i)}.
} 
\\
Case {\em(ii)}:
$\VH = \VS \ordSum \VF$,
$\{x, y\}\cap S \not= \emptyset$,
$\{x, y\}\cap F \Diff \{0\} \not= \emptyset$:
This leads to 2 subcases.
These are verified using elementary properties of $\VS \ordSum \VF$ as follows:
\\
Subcase {\em(ii)(a)}: $x \in S, y \in F \Diff \{0\}$:
we have:
\[
i(x) \rImp i(y) \rImp i(x \rImp y) = i(x) \rImp i(y) \rImp i(y) = i(x) \rImp 0 = 0.
\]
Subcase {\em(ii)(b)}: $x \in F \Diff \{0\}, y \in S$:
we have:
\[
i(x) \rImp i(y) \rImp i(x \rImp y) = i(x) \rImp i(y) \rImp 0 = 0.
\]
This completes case {\em(ii)}.
\Done
\section{Double Negation Translations}\label{sec:neg-trans}

In this section we undertake an algebraic study of the syntactic
translations known as double negation translations (or negative translations).
Throughout this section all our pocrims will be bounded.
We will view $1$ as a constant in the signature for bounded pocrims and so a
homomorphism $f$ must satisfy $f(1) = 1$. In the semantics, $1$ will always be
interpreted as $1$, so assignments will be functions with domain $\Var$ rather
than $\Var \cup\{1\}$.

\subsection{The Double Negation Mapping}

If $\VP$ is a pocrim, let $N = \Im(\Not) = \{\Not x \ST x \in P\}$.
Since $\delta(\Not x) = \Not x$, $\Im(\delta) = N$.
Clearly $\{0, 1\} \subseteq N$ and $N$ is closed under $\rImp$,
since $\Not x \rImp \Not y = \Not(\Not x + y)$.
In general, $N$ is not closed under addition and hence is not a subpocrim
and $\delta$ does not respect either $+$ or $\rImp$:

\begin{Example}
There is a pocrim $\VU$ with elements $0 < a < b < c < 1$
and with $+$, $\rImp$ and $\delta$ as follows:
\[
\begin{array}{l@{\quad\quad}l@{\quad\quad}l}
\begin{array}{c|ccccc}
   {+} & 0 & a & b & c & 1 \\\hline
    0  & 0 & a & b & c & 1 \\
    a  & a & b & b & 1 & 1\\
    b  & b & b & b & 1 & 1\\
    c  & c & 1 & 1 & 1 & 1\\
    1  & 1 & 1 & 1 & 1 & 1
\end{array}
&
\begin{array}{c|ccccc}
   {\rImp} & 0 & a & b & c & 1 \\\hline
    0      & 0 & a & b & c & 1 \\
    a      & 0 & 0 & a & c & c \\
    b      & 0 & 0 & 0 & c & c \\
    c      & 0 & 0 & 0 & 0 & a \\
    1      & 0 & 0 & 0 & 0 & 0
\end{array}
&
\begin{array}{c|c}
   \multicolumn{2}{c}{\delta}  \\
   \hline
    0      & 0 \\
    a      & a \\
    b      & a \\
    c      & c \\
    1      & 1
\end{array}
\end{array}
\]
So, in $\VU$,  $\delta(a \rImp b) = a \not= 0 = \delta(a) \rImp \delta(b)$,
$\delta(a + a) = a \not= b = \delta(a) + \delta(a)$
and $\delta(\delta(a) + \delta(a)) \not= \delta(a) + \delta(a)$.
The image of negation is
$\{0, a, c, 1\}$, which is not closed under addition, since $a + a = b$.
\end{Example}

The situation in a hoop is much more satisfactory. To describe it,
we first make the following definition:

\begin{Definition}\label{def:cl-cl}
If $\VH$ is a bounded hoop, the {\em involutive replica}, $\ClCl{\VH}$,
of $\VH$ is $\VH/\theta$, where $\theta$ is the smallest congruence
such that $x \mathrel{\theta} \delta(x)$ for all $x \in H$.
\end{Definition}

$\VH \mapsto \ClCl{\VH}$ is the objects part of a functor from the category
of bounded hoops to the category of involutive hoops and every
homomorphism from $\VH$ to an involutive hoop factors
uniquely through $\ClCl{\VH}$.

\begin{Theorem}\label{thm:double-neg}
If $\VH$ is a bounded hoop, then the double negation mapping,
$\delta$, is a homomorphism $\VH \To \VH$.
Moreover, if
$p : \VH \To \ClCl{\VH}$ is the natural projection,
then $p$ factors as
$p = i \circ \delta$ where $i : \Im(\delta) \To \ClCl{\VH}$
is an isomorphism.
\end{Theorem}
\Proof
By Theorem~\ref{thm:de-morgan}, we have:
\begin{align*}
&\delta(x) + \delta(y) = \Not(x \rImp \Not y) = \delta(x + y) \\
&\delta(x) \rImp \delta(y) = \Not(\delta(x) + \Not y) = \delta(x) \rImp \delta(y)
\end{align*}
As $\delta$ fixes the constants $0$ and $1$, this proves that $\delta$ is a
homomorphism.  The claim about $p$ is equivalent to the claim that
$\Ker(\delta) = \Ker(p)$.  Now $\ClCl{\VH}$ is the quotient $\VH/\theta$
where $\theta$ is the smallest congruence such that $x \mathrel{\theta}
\delta(x)$ for all $x \in H$.  As $\delta$ is an idempotent endomorphism, it is
not difficult to see that $x \mathrel{\theta} y$ iff $\delta(x) = \delta(y)$.
So we have $\Ker(\delta) = \{x \ST \delta(x) = 0 \} = \{x \ST \delta(x) =
\delta(0)\} = \{x \ST x \mathrel{\theta} 0\} = \Ker(p)$.
\Done

\subsection{Semantics for Double Negation Translations}

Beginning with Kolmogorov \cite{Kolmogorov(25)}, logicians have
studied {\em double negation translations} that represent classical logic in
intuitionistic logic. Kolmogorov's translation inductively replaces
every subformula of a formula by its double negation. Subsequent authors have
devised more economical translations:
Gentzen's translation \cite{Gentzen(33)} applies double negation
to atomic formulas only and Glivenko's translation \cite{Glivenko(29)}
just applies double negation to a formula without changing its
internal structure.

We wish to undertake an algebraic analysis of translations such as
the various double negation translations. We will view the translations
as variant semantics and so we need a framework to compare
different semantics.

\begin{Definition}
Let $\BPoc$ be the category of bounded pocrims and homomorphisms
and let $\Set$ be the category of sets.
Given any set $X$, let $H_X : \BPoc \To \Set$
be the functor that maps a pocrim $\VP$
to $\Hom_{\Set}(X, P)$, i.e.,
the set of all functions from $X$ to $P$, and 
maps a homomorphism $h : \VP \To \VQ$ to
$f \mapsto h \circ f : \Hom_{\Set}(X, P) \To \Hom_{\Set}(X, Q)$.
Now let $\Ass = H_{\Var}$ and $\Sem = H_{\lL}$.
We define a {\em semantics} to be a natural transformation
$\mu : \Ass \To \Sem$.
\end{Definition}

So given a bounded pocrim $\VP$, $\Ass(\VP)$ denotes the set of assignments
$\alpha : \Var \To P$ , while $\Sem(\VP)$ denotes the set of all possible
functions $s : {\lL} \To P$.
A semantics $\mu$ is a family of functions $\mu_{\VP}$
indexed by bounded pocrims $\VP$ such that $\mu_{\VP} : \Ass(\VP) \To \Sem(\VP)$
and such that for any homomorphism $f : \VP \To \VQ$ the following diagram commutes.
\[
\begin{CD}
\Ass(\VP) @>\Ass(f)>> \Ass(\VQ) \\
@VV\mu_{\VP}V          @VV\mu_{\VQ}V\\
\Sem(\VP) @>\Sem(f)>> \Sem(\VQ) \\
\end{CD}
\]

The standard semantics $\sSem$ is the one used to define bounded
validity in Section~\ref{sec:alg-sem}: it simply uses the given assignment
$\alpha:\Var \To P$ to give values to the variables in a formula in $\lL$
and then calculates its value interpreting $1$, $\iAnd$ and $\Lolly$
as $1$, $+$ and $\rImp$ respectively:
\begin{align*}
\sSem_{\VP}(\alpha)(V_i) &= \alpha(V_i) \\
\sSem_{\VP}(\alpha)(1) &= 1 \\
\sSem_{\VP}(\alpha)(A \iAnd B) &= \sSem_{\VP}(\alpha)(A) + \sSem_{\VP}(\alpha)(B) \\
\sSem_{\VP}(\alpha)(A \Lolly B) &= \sSem_{\VP}(\alpha)(A) \rImp \sSem_{\VP}(\alpha)(B)
\end{align*}
The Kolmogorov translation corresponds to a semantics $\kSem$ defined
like $\sSem$, but applying double negation to everything in sight:
\begin{align*}
\kSem_{\VP}(\alpha)(V_i) &= \delta(\alpha(V_i)) \\
\kSem_{\VP}(\alpha)(1) &= 1 \\
\kSem_{\VP}(\alpha)(A \iAnd B)) &=
   \delta(\kSem_{\VP}(\alpha)(A) + \kSem_{\VP}(\alpha)(B)) \\
\kSem_{\VP}(\alpha)(A \Lolly B)) &=
   \delta(\kSem_{\VP}(\alpha)(A) \rImp \kSem_{\VP}(\alpha)(B))
\end{align*}
It is easily verified that $\sSem$ and $\kSem$ are indeed natural
transformations $\Ass \To \Sem$.
The Gentzen and Glivenko translations correspond to semantics obtained
by composing the standard semantics with double negation:
\begin{align*}
\ggSem &= \sSem \circ \delta^{\Var} \\
\glSem &= \delta^{\lL} \circ \sSem
\end{align*}
where $\delta^X$ denotes the natural transformation from
$H_X = \Hom_{\Set}(X, \cdot)$ to itself with $\delta^X_{\VP} = f \mapsto \delta \circ f$.
It is clear from Theorem~\ref{thm:double-neg} that the Kolmogorov, Gentzen
and Glivenko semantics all agree when restricted to hoops.

\begin{Definition}
Let $\cC$ be a class of bounded pocrims, we say that a semantics $\mu$
is a {\em double negation semantics} for $\cC$ if the following
conditions hold:
\begin{description}
\item[(DNS1)] If $\VP \in \cC$ is involutive, then $\mu_{\VP} = \sSem_{\VP}$.
\item[(DNS2)] Given a formula $A$, if, for every involutive $\VP \in \cC$
and every $\alpha : \Var \To P$, we have:
\[
\sSem_{\VP}(\alpha)(A) = 0,
\]
then, for every $\VP \in \cC$ and every $\alpha : \Var \To P$, we have:
\[
\mu_{\VP}(\alpha)(A) = 0.
\]
\item[(DNS3)] $\delta^{\lL} \circ \mu = \mu$.
\end{description}
\end{Definition}
\begin{Remark}\label{rmk:dns}
Let us write $\Th(\cC)$ for the {\em theory} of a class of pocrims, i.e., the
set of all formulas $A$ such that $\sSem_{\VP}(\alpha)(A) = 0$ for every
$\alpha : \Var \To P$ where $\VP \in \cC$. The above definition can be seen to
agree with the usual syntactic definition of a double negation translation due
to Troelstra \cite{Troelstra(73)}, provided $\Th(\cI) = \Th(\cC) + \axDNE$,
where $\cI$ comprises the involutive pocrims in $\cC$.  See
\cite{arthan-oliva14a} for more information about the various syntactic double
negation translations in $\ALm$ and its extensions.
\end{Remark}

\begin{Theorem}
The Kolmogorov semantics,
$\kSem$, the Gentzen semantics, $\ggSem$, and the Glivenko semantics,
$\glSem$, are double negation semantics for any class of hoops.
\end{Theorem}
\Proof
(DNS1) and (DNS3) are clear for $\mu = \glSem = \delta^{\lL} \circ \sSem$,
since $\delta^{\lL}_{\VH} = \Id(H)$ when $\VH$ is involutive
and $\delta^{\lL} \circ \delta^{\lL} = \delta^{\lL}$.
Also (DNS2) holds for $\mu = \ggSem$ in the class of hoops, since, if $\VH$ is a hoop, then $\Im(\delta)$ is an involutive subhoop, and
for any $\alpha : \Var \To \VH$, we have:
\[
\ggSem_{\VH}(\alpha) = (\sSem_{\VH} \circ \delta^{\Var})(\alpha) = \sSem_{\VH}(\delta \circ \alpha)
  = \sSem_{\Im(\delta)}(\delta \circ \alpha).
\]
Now it is easy to see using Theorem~\ref{thm:double-neg}, that if $\VH$ is a hoop, then we have:
\[
\kSem_{\VH} = \ggSem_{\VH} = \glSem_{\VH}
\]
Hence (DNS1), (DNS2) and (DNS3) hold for any of the three translations in any class of hoops.
\Done

\begin{Lemma}\label{lma:dn-imp}
Any pocrim satisfies $\delta(a \rImp b) \ge \delta(a) \rImp \delta(b)$.
\end{Lemma}
\Proof
It is easy to see that
$(*)$ if $x + y = 1$, then $x \ge \Not y$ and
$(**)$ if $x \ge y$, then $\Not x + y = 1$.
Combining $(*)$ and $(**)$, we have $(***)$ if $x + y = 1$ then $\delta(x) + y = 1$.
Hence:
\begin{align*}
a + (a \rImp b) &\ge b \\
a + \Not b + (a \rImp b) &= 1 \tag*{$(**)$} \\
\delta(a) + \Not b + \delta(a \rImp b) &= 1 \tag*{$2 \times \mbox{$(***)$}$} \\
\delta(a) + \delta(a \rImp b) &\ge \delta(b) \tag*{$(*)$} \\
\delta(a \rImp b) &\ge \delta(a) \rImp \delta(b)  \tag*{\Done}
\end{align*}

\begin{Theorem}
The Kolmogorov semantics,
$\kSem$, is a double negation semantics for $\ALi$.
\end{Theorem}
\Proof
(DNS1) and (DNS3) are easy to verify.
For (DNS2), by Theorem~\ref{thm:sound-complete} it is enough to prove
that, if $\ALc$ proves $A$, then, for any pocrim $\VP$ and any
$\alpha : \Var \To P$, $\kSem_{\VP}(\alpha)(A) = 0$.
So let $\VP$ and $\alpha : \Var \To P$, be given.
We prove this by induction on a proof of $A$.
The induction has an inductive step corresponding to our single inference
rule and a base case for each of the axiom schemata used to define $\ALc$. \\[1mm]
{\em Modus ponens}: by the inductive hypothesis,
we are given that $\ALc$ proves $B$ and $B \Lolly A$.
Let $a = \kSem_{\VP}(\alpha)(A)$ and $b = \kSem_{\VP}(\alpha)(B)$
and note that from the definition of $\kSem$ this means $a \in \Im(\delta)$.
We want to show that $a = 0$. By the inductive hypothesis
$b = 0$ and $\delta(b \rImp a) = 0$, but then as $a \in \Im(\delta)$
and using Lemma~\ref{lma:dn-imp},
we have $a = \delta(a) = \delta(b) \rImp \delta(a) = 0$. \\[1mm]
For the axiom schemata, we have to show that if $A$ is an instance
of one of the schemata, then the semantic value $X = \kSem_{\VP}(\alpha)(A)$
of the instance is equal to 0.
We will make frequent and tacit use of the facts that
$x + y \rImp z = x \rImp y \rImp z$, $x \ge \delta(x)$, that
$x \rImp y \ge \delta(x) \rImp \delta(y)$ and that, if $x \in \Im(\delta)$,
then $x = \delta(x)$. \\[1mm]
$\axComp$: In this case, $X = \delta(\delta(a \rImp b) \rImp \delta(\delta(b \rImp c) \rImp \delta(a \rImp c))$, for some $a$, $b$ and $c$, and we have:
\begin{align*}
\delta(\delta(a \rImp b) \rImp \delta(\delta(b \rImp c) \rImp \delta(a \rImp c)) &\le
  \delta((a \rImp b) \rImp \delta(b \rImp c) \rImp \delta(a \rImp c)) \\
   &\le \delta((a \rImp b) \rImp (b \rImp c) \rImp (a \rImp c)) \\
   &\le \delta(0) = 0
\end{align*} \\[1mm]
$\axComm$: $X = \delta(\delta(a+b) \rImp \delta(b+a))$, for some $a$ and $b$,
so $X = \delta(\delta(a+b) \rImp \delta(a+b)) = \delta(0) = 0$. \\[1mm]
$\axCurry$: $X = \delta(Y \rImp Z)$ where
$Y = \delta(\delta(a + b) \rImp c)$ and $Z = \delta(a \rImp \delta(b \rImp c))$,
for some $a$, $b$ and $c$, and it is enough to prove $Y \ge Z$. We have:
\begin{align*}
\delta(\delta(a + b) \rImp c)
   &\ge \delta(a + b \rImp c)
      \tag*{(as $a + b \ge \delta(a + b)$)} \\
   &= \delta(a \rImp b \rImp c)
      \\
   &\ge \delta(a \rImp \delta(b \rImp c))
      \tag*{(as $b \rImp c \ge \delta(b \rImp c)$.}
\end{align*} \\[1mm]
$\axUncurry$: $X = \delta(Y \rImp Z)$ where
$Y = \delta(a \rImp \delta(b \rImp c))$ and
$Z = \delta(\delta(a + b) \rImp c))$, for some $a, b, c \in \Im(\delta)$,
and it is enough to prove $Y \ge Z$. We have:
\begin{align*}
\delta(a \rImp \delta(b \rImp c)) &\ge
   \delta(a \rImp \delta(b) \rImp \delta(c))
      \tag*{(Lemma~\ref{lma:dn-imp})} \\
   &= \delta(a \rImp b \rImp c)
      \tag*{(as $b, c \in \Im(\delta)$)} \\
   &= \delta(a + b \rImp c) \\
   &\ge \delta(a + b) \rImp \delta(c)
      \tag*{(Lemma~\ref{lma:dn-imp})} \\
   &= \delta(a + b) \rImp c
      \tag*{(as $c \in \Im(\delta)$)} \\
   &\ge \delta(\delta(a + b) \rImp c).
\end{align*} \\[1mm]
$\axWk$: $X = \delta(\delta(a + b) \rImp a)$ where
$a \in \Im(\delta)$. We have:
\begin{align*}
\delta(\delta(a + b) \rImp a)
	& = \delta(\delta(a + b) \rImp \delta(a)) \\
	& \le \delta((a + b) \rImp a) = \delta(0) = 0 
\end{align*} \\
$\axEFQ$: For some $a$, $X = \delta(1 \rImp a)$, so $X = \delta(0) = 0$. \\[1mm]
$\axDNE$: For some $a \in \Im(\delta)$, $X = \delta(\delta(a) \rImp a) = \delta(a \rImp a) = \delta(0) = 0$.
\Done

\begin{Example}\label{eg:q-six}
Consider the pocrim $\VQ_6$ with six elements $0 < p < q < r < s < 1$
and with $+$, $\rImp$ and $\delta$ as shown in the following tables:
\[
\begin{array}{l@{\quad\quad}l@{\quad\quad}l}
\begin{array}{c|cc|cc|c|c}
	+ & 0 & p & q & r & s & 1 \\
	\hline
	0 & 0 & p & q & r & s & 1 \\
	p & p & p & r & r & s & 1 \\
	\hline
	q & q & r & r & r & 1 & 1 \\
	r & r & r & r & r & 1 & 1 \\
	\hline
	s & s & s & 1 & 1 & 1 & 1 \\
	\hline
	1 & 1 & 1 & 1 & 1 & 1 & 1
\end{array}
&
\begin{array}{c|cc|cc|c|c}
	\rImp & 0 & p & q & r & s & 1 \\
	\hline
	0 & 0 & p & q & r & s & 1 \\
	p & 0 & 0 & q & q & s & 1 \\
	\hline
	q & 0 & 0 & 0 & p & s & s \\
	r & 0 & 0 & 0 & 0 & s & s \\
	\hline
	s & 0 & 0 & 0 & 0 & 0 & q \\
	\hline
	1 & 0 & 0 & 0 & 0 & 0 & 0
\end{array}
&
\begin{array}{c|c}
	\multicolumn{2}{c}{\delta} \\
	\hline
	0 & 0 \\
	p & 0 \\
	\hline
	q & q \\
	r & q \\
	\hline
	s & s \\
	\hline
	1 & 0
\end{array}
\end{array}
\]
$\VQ_6$ is not involutive, as $\delta(x) = x$ fails for $x \in \{p, r\}$.
In $\VQ_6$, double negation is an implicative homomorphism:
$\Not\Not x \rImp \Not\Not y = \Not\Not(x \rImp y)$ for all $x, y$.
Double negation is not quite an additive homomorphism in $\VQ_6$:
$\Not\Not x + \Not\Not y = \Not\Not(x + y)$
unless $\{x, y\} \subseteq \{q, r\}$, in which case
$\Not\Not x + \Not\Not y = r > q = \Not\Not(x + y)$.
As indicated by the block decomposition of the operation tables,
there is a homomorphism $h : \VQ_6 \To \VQ_4$, where
$\VQ_4$ is as discussed in Example~\ref{eg:q-four}.
The kernel congruence of $h$ has equivalence classes
$\{0, p\}$, $\{q, r\}$, $\{s\}$ and $\{1\}$ which are mapped by $h$
to $0$, $u$, $v$, $1$ respectively in $\VQ_4$.
\end{Example}

\begin{Theorem}\label{thm:gentzen-glivenko-not-ali}
(i) The Gentzen semantics $\ggSem$ is not a double negation semantics for any class of
pocrims that contains the pocrim $\VQ_6$ of Example~\ref{eg:q-six}. (ii) The Glivenko semantics $\glSem$ is not a double negation semantics for any class
of pocrims that contains the pocrim $\VP_4$ of Example~\ref{eg:p-four}.
\end{Theorem}
\Proof {\em(i):} We show that (DNS2) does not hold for $\ggSem$ in $\VQ_6$.
Let $V, W \in \Var$ and let $A$ be the formula
 $ (V \iAnd W)\Lnot\Lnot \Lolly (V \iAnd W)$.
$A$ is an instance of $\axDNE$ and so, by Theorem~\ref{thm:sound-complete},
$\sSem_{\VP}(\alpha)(A) = 0$, for any involutive pocrim $\VP$ and any
$\alpha : \Var \To P$. Thus (DNS2) requires $\ggSem_{\VQ_6}(\alpha)(A) = 0$
for any $\alpha:\Var \To Q_6$.
However, if $\alpha(V) = \alpha(W) = r$, we have:
\begin{align*}
\ggSem_{\VQ_6}(\alpha)(A) &=
   \delta(\delta(r) + \delta(r)) \rImp \delta(r) + \delta(r) \\
   & = \delta(q + q) \rImp q + q  \\
   &= \delta(r) \rImp r =
   q \rImp r = s \not = 0.
\end{align*}
\\
{\em(ii):}
we argue as in the proof of (A), but taking $A$ to be
$V\Lnot\Lnot \Lolly V$. Then, if $\alpha(V) = q$, we have:
\begin{align*}
\glSem_{\VP}(\alpha)(A) &= \delta(\delta(q) \rImp q) \\
   &= \delta(p\rImp q) = \delta(p) = p \not= 0. \tag*{\Done}
\end{align*}

\begin{Theorem}\label{thm:gliv-gentz}
Let $\cC_1$ comprise the two pocrims $\VP_4$ and $\VL_3$
of Examples~\ref{eg:vl} and~\ref{eg:p-four} and
let $\cC_2$ comprise the two pocrims $\VQ_6$ and $\VQ_4$
of Examples~\ref{eg:q-six} and~\ref{eg:q-four}.
Then:
\\
{\em(i)} The Gentzen semantics, $\ggSem$, is a double negation semantics for $\cC_1$, but the Glivenko semantics, $\glSem$, is not.
\\
{\em(ii)} The Glivenko semantics, $\glSem$, is a double negation semantics for $\cC_2$, but the Gentzen semantics, $\ggSem$, is not.
\end{Theorem}
\Proof
{\em(i)}:
By Theorem~\ref{thm:gentzen-glivenko-not-ali}, $\glSem$ is not a double negation semantics for $\cC_1$.
As for $\ggSem$, (DNS1) is easily verified.
For (DNS3) and (DNS2), note that for any $\alpha : \Var \To P_4$, we
have:
\[
\ggSem_{\VP_4}(\alpha) = (\sSem_{\VP_4} \circ \delta^{\Var})(\alpha) = \sSem_{\VP_4}(\delta \circ \alpha) = \sSem_{\VL_3}(\delta \circ \alpha)
\]
where in the last expression we have identified $\VL_3$
with the subpocrim of $\VP_4$ whose universe is $\Im(\delta)$.
Thus evaluation under $\ggSem$ with an assignment in any pocrim in $\cC_1$ is equivalent
to evaluation under the standard semantics, $\sSem$,
with an assignment in the involutive pocrim $\VL_3$. (DNS3) and (DNS2) follow immediately from this.
\\
{\em(ii)}:
By Theorem~\ref{thm:gentzen-glivenko-not-ali}, $\ggSem$ is not a double negation semantics for $\cC_2$.
As for $\glSem$, (DNS1) and (DNS3) are immediate from the definition of $\glSem$.
For (DNS2), let $A$ be a formula, such that $\sSem_{\VQ_4}(\alpha)(A) = 0$,
for any assignment $\alpha : \Var \To \VQ_4$. As $\VQ_4$ is the only
involutive pocrim in $\cC_2$, we must show that $\glSem_{\VP}(\alpha)(A) = 0$ for $\VP \in \cC_2$ under any assignment $\alpha : \Var \To P$.
This is easy to see for $\VP = \VQ_4$, since the Glivenko semantics is
the double negation of the standard semantics and $\VQ_4$ is involutive.
As for $\VP = \VQ_6$, let $\alpha : \Var \To \VQ_6$ be given.
As discussed in Example~\ref{eg:q-six}, there is a homomorphism
$h : \VQ_6 \To \VQ_4$, so, as $\sSem$ is a natural transformation,
the following diagram commutes:
\[
\begin{CD}
\Ass(\VQ_6) @>\Ass(h)>> \Ass(\VQ_4) \\
@VV\sSem_{\VQ_6}V          @VV\sSem_{\VQ_4}V\\
\Sem(\VQ_6) @>\Sem(h)>> \Sem(\VQ_4) \\
\end{CD}
\]
Hence, by the assumption on $A$,  we have:
\[
(h \circ \sSem_{\VQ_6}(\alpha))(A) = \sSem_{\VQ_4}(h \circ \alpha)(A) = 0
\]
So $\sSem_{\VQ_6}(\alpha)(A) \in h^{-1}(0) = \{0, p\}$. As
$\delta(0) = \delta(p) = 0$, we can conclude:
\begin{align}
\glSem_{\VQ_6}(\alpha)(A) = \delta(\sSem_{\VQ_6}(\alpha)(A)) = 0. \tag*{\Done}
\end{align}

\begin{Remark}
Taken with the following lemma and Remark~\ref{rmk:dns},
Theorem~\ref{thm:gliv-gentz} implies the existence of logics extending $\ALi$
in which the syntactic Gentzen translation meets Troelstra's requirements on a
double negation translation but the syntactic Glivenko translation does not and vice
versa.
\end{Remark}
\begin{Lemma} If $\VL_3$, $\VP_4$, $\VQ_4$ and $\VQ_6$ are as in Theorem~\ref{thm:gliv-gentz}, then:
\begin{align*}
\Th(\VL_3) &= \Th(\VP_4) + \axDNE \\
\Th(\VQ_4) &= \Th(\VQ_6) + \axDNE.
\end{align*}
\end{Lemma}
\Proof
For the first equation, the
right-to-left inclusion holds because identities are preserved
in subalgebras.
For left-to-right, let us write $\VP \models A$ to mean $\sSem_{\VP}(\alpha)(A)
= 0$ for every $\alpha : \Var \To P$.  Assume $\VL_3 \models A$ and let $W_1,
\ldots, W_k$ be the variables occurring in $A$. Define $B$ to be
$(W_1\Lnot\Lnot \Lolly W_1) + \ldots + (W_k\Lnot\Lnot \Lolly W_k)$. We claim that
$\sSem_{\VP_4}(\alpha)(B \iAnd B \Lolly A) = 0$ for every $\alpha :
\Var \To P_4$, so that as $\ALm + \axDNE$ proves $B$, $\Th(\VP_4) + \axDNE$
proves $A$. To see this let an assignment $\alpha : \Var \To P_4$ be given.
Then either {\em(i)} $\Im(\alpha) \subseteq \{0,  p, 1\}$, in which case
$\sSem_{\VP_4}(\alpha)(A) = 0$, since $\alpha$ is an assignment into a subpocrim isomorphic
to $\VL_3$ and $\VL_3 \models A$ by assumption, or {\em(ii)} $\alpha(W_i) = q$
for some $i$, but then $\sSem_{\VP_4}(\alpha)(W_i\Lnot\Lnot \Lolly W_i) = q$
and so $\sSem_{\VP_4}(\alpha)(B \iAnd B) \ge q + q = 1$.  In both cases, we
have that $\sSem_{\VP_4}(\alpha)(B \iAnd B \Lolly A) = 0$,
proving the claim. The proof of the second equation is similar using the facts that identities are
preserved in quotient algebras and that, if $\VQ_4 \models A$ and
$\alpha : \Var \To Q_6$, then $\sSem_{\VQ_6}(\alpha)(A) \in \{0, p\}$, implying
that $\VQ_6 \models (A\Lnot\Lnot \Lolly A) \Lolly A$.
\Done

\section{Concluding Remarks}
The axiom $A \iAnd (A \Lolly B) \Lolly B \iAnd (B \Lolly A)$, that we call
$\axCWC$, characterizes what seems to us to be an important landmark between
affine logic, in which using an assumption destroys it, and standard logic, in
which we may use an assumption as many times as we please.
The importance of this axiom is reflected algebraically in the rich properties
enjoyed by hoops, the algebraic models of $\axCWC$, when
compared with the algebraic models of general affine logic, namely pocrims.
However, many of these properties depend on algebraic laws whose
derivations involve extremely intricate applications of the identity $x + (x
\rImp y) = y + (y \rImp x)$ that corresponds to the axiom $\axCWC$. The methods
of the present paper mitigate the problem of finding these derivations in many
cases of interest.

Our original interest in {\L}ukasiewicz logic arose from work on {\em continuous
logic} \cite{ben-yaacov-pedersen09}.  In \cite{arthan-oliva12} we investigate
the natural algebraic models for continuous logic and an intuitionistic
analogue. These models comprise specialisations of hoops that we call {\em
coops} which admit a halving operator $x \mapsto x/2$ satisfying the law $x/2 =
x/2 \rImp x$. There is a characterization of subdirectly irreducible coops
very like Blok and Ferreirim's result for hoops and
the  method for proving identities of the present paper
carries over straightforwardly.
Using it, one may prove, for example, the following ``De Morgan'' identity:
$\Not(x/2) = 1/2 + (\Not x)/2$.

Bova and Montagna have shown that the quasi-equational theory of
commutative GBL-algebras is PSPACE-complete and conjecture that
the equational theory is also PSPACE-complete. Our indirect method
of proof provides a heuristic that proves to be very successful on
simple formulas with just a few variables. One might conjecture
that the decision problem for a fixed number of variables admits
a more tractable decision procedure.
A problem would be to give a tractable description
of the structure of the free hoop on $n$ generators.  Berman and Blok
\cite{berman-blok04} have studied free $k$-potent hoops, but the assumption of
$k$-potency is quite a strong one: e.g., a coop $\VC$ is $k$-potent iff $C =
\{0\}$.

Semantic methods of some sort
are the only way of obtaining results such as Theorem~\ref{thm:gliv-gentz} that
delimit the applicability of a given syntactic translation.  Hyland
\cite{Hyland02} gives a semantic account of double negation translations in
categorical terms.
It would be interesting to attempt to integrate the
categorical approach with the algebraic approach of the present paper. \\[4mm]
{\bf Acknowledgments}
We thank: Franco Montagna for drawing our attention to \cite{Bova:2009}
and for a manuscript proof that the equational theory of
commutative GBL-algebras is a conservative extension of that of hoops;
George Metcalfe for encouraging remarks and for pointers to the literature;
and Isabel Ferreirim for helpful correspondence about the theory of
hoops.

\bibliographystyle{plain}

\bibliography{../references}

%
%

\end{document}